\documentclass[a4paper,12pt]{amsart}
\usepackage{amsfonts}
\usepackage{mathrsfs}
\usepackage{amssymb}
\usepackage{ifthen}
\usepackage{graphicx}
\usepackage{tablists}
\nonstopmode \numberwithin{equation}{section}
\usepackage[usenames]{color}
\usepackage{enumerate}
\usepackage{tikz}
\newtheorem{thm}[equation]{Theorem}
\newtheorem{lem}[equation]{Lemma}
\newtheorem{cor}[equation]{Corollary}

\newtheorem{cl}{Claim}[section]
\newtheorem{ca}{Case}[section]
\newtheorem{sca}{Subcase}[section]
\newtheorem{scl}[section]{Subclaim}
\newtheorem{conj}[equation]{Conjecture}

\theoremstyle{definition}
\newtheorem{defn}[equation]{Definition}
\newtheorem{prob}[equation]{Problem}
\newtheorem{op}[equation]{Open Problem}
\newtheorem{ques}[equation]{Question}
\newtheorem{rem}[equation]{Remark}
\newtheorem{exam}[equation]{Example}

\newcounter {own}
\def\theown {\thesection       .\arabic{own}}

\makeatletter
\@namedef{subjclassname@2020}{\textup{2020} Mathematics Subject Classification}
\makeatother

\newenvironment{pf}[1][]{%
 \vskip 3mm
 \noindent
 \ifthenelse{\equal{#1}{}}%
  {{\slshape Proof. }}%
  {{\slshape #1.} }%
 }%
{\qed\bigskip}

\newcounter{alphabet}


\newcommand{\IR}{{\mathbb R}}

\newcommand{\id}{{\operatorname{id}}}

\newcommand{\diam}{{\operatorname{diam}}}

\newcommand{\dist}{{\operatorname{dist}}}

\def\be{\begin{equation}}
\def\ee{\end{equation}}

\newcommand{\ben}{\begin{enumerate}}
\newcommand{\een}{\end{enumerate}}

\newcommand{\blem}{\begin{lem}}
\newcommand{\elem}{\end{lem}}
\newcommand{\bthm}{\begin{thm}}
\newcommand{\ethm}{\end{thm}}
\newcommand{\bcor}{\begin{cor}}
\newcommand{\ecor}{\end{cor}}
\newcommand{\beg}{\begin{exam}}
\newcommand{\eeg}{\end{exam}}
\newcommand{\begs}{\begin{examples}}
\newcommand{\eegs}{\end{examples}}
\newcommand{\bdefe}{\begin{defn}}
\newcommand{\edefe}{\end{defn}}
\newcommand{\bprob}{\begin{prob}}
\newcommand{\eprob}{\end{prob}}
\newcommand{\bques}{\begin{ques}}
\newcommand{\eques}{\end{ques}}
\newcommand{\bei}{\begin{itemize}}
\newcommand{\eei}{\end{itemize}}
\newcommand{\bcon}{\begin{conj}}
\newcommand{\econ}{\end{conj}}
\newcommand{\bop}{\begin{op}}
\newcommand{\eop}{\end{op}}

\newcommand{\bca}{\begin{ca}}
\newcommand{\eca}{\end{ca}}
\newcommand{\bsca}{\begin{sca}}
\newcommand{\esca}{\end{sca}}

\newcommand{\bcl}{\begin{cl}}
\newcommand{\ecl}{\end{cl}}

\newcommand{\bscl}{\begin{scl}}
\newcommand{\escl}{\end{scl}}

\newcommand{\bcons}{\begin{conjs}}
\newcommand{\econs}{\end{conjs}}
\newcommand{\bprop}{\begin{propo}}
\newcommand{\eprop}{\end{propo}}
\newcommand{\br}{\begin{rem}}
\newcommand{\er}{\end{rem}}
\newcommand{\brs}{\begin{rems}}
\newcommand{\ers}{\end{rems}}
\newcommand{\bo}{\begin{obser}}
\newcommand{\eo}{\end{obser}}
\newcommand{\bos}{\begin{obsers}}
\newcommand{\eos}{\end{obsers}}
\newcommand{\bpf}{\begin{pf}}
\newcommand{\epf}{\end{pf}}
\newcommand{\ba}{\begin{array}}
\newcommand{\ea}{\end{array}}
\newcommand{\beq}{\begin{eqnarray}}
\newcommand{\beqq}{\begin{eqnarray*}}
\newcommand{\eeq}{\end{eqnarray}}
\newcommand{\eeqq}{\end{eqnarray*}}

\newcounter{minutes}
\divide\time by 60
\newcounter{hours}
\multiply\time by 60 \addtocounter{minutes}{-\time}

\begin{document}

\bibliographystyle{amsplain}
\title{Gromov hyperbolicity and uniformity}

\author{Qingshan Zhou 
}
\address{Qingshan Zhou, School of Mathematics, Foshan University,  Foshan, Guangdong 528000, People's Republic
of China} \email{qszhou1989@163.com; q476308142@qq.com}

\def\thefootnote{}
\footnotetext{ \texttt{\tiny File:~\jobname .tex,
          printed: \number\year-\number\month-\number\day,
          \thehours.\ifnum\theminutes<10{0}\fi\theminutes}
} \makeatletter\def\thefootnote{\@arabic\c@footnote}\makeatother

\date{}

\begin{abstract} In this paper, we first exhibit an example which answers negatively to a question posed by Herron et al. in 2008. That is, there is a rectifiably connected open set in a proper geodesic metric space so that it is Gromov hyperbolic with respect to the quasihyperbolic metric, and its Gromov boundary endowed with a visual metric is naturally bilipschitz equivalent to its metric boundary, but it is not uniform.

Moreover, two positive answers are given under additional hypotheses. It is shown that a metric space is uniform if and only if it is Gromov hyperbolic and the natural map between its Gromov closure and its metric completion exists and satisfies a relatively quasisymmetric condition. Both bounded and unbounded cases are in our considerations via two different uniformization theorems developed recently by Bonk et al. and the first author, respectively. Finally, we characterize uniform metric spaces with uniformly perfect boundaries in terms of Gromov hyperbolicity and Teichm\"uller's displacement property.
\end{abstract}

\thanks{Qingshan Zhou was supported by the National Natural
Science Foundation of China (No. 12671090).}

\maketitle{} \pagestyle{myheadings} \markboth{Qingshan Zhou} {Gromov hyperbolicity and uniformity}

\section{Introduction and main results}\label{sec-1}
\subsection{Gromov hyperbolicity}
It is known that hyperbolic $n$-space $\mathbb{H}^n$ is a complete simply connected Riemannian manifold with all sectional curvatures $-1$ at every point. Historically, the main motivation for the theory of negatively curved metric spaces came from geometric group theory, in particular, in the study of groups acting on manifolds of negative curvature. In \cite{Gr87}, Gromov attempted to give a geometric characterization of these groups in terms of their Cayley graphs, and realized that the asymptotic properties of $\mathbb{H}^n$ can be ensured by using a simple condition for quadruples of points, which is now used to define Gromov hyperbolic spaces.

The class of Gromov hyperbolic spaces forms a large and much studied class of metric spaces and has influenced much research, see \cite{BBGS17,BrHa,BDS,Ca,FS11,ZLR}. For example, it was shown by Bonk and Schramm in \cite{BS} that a $\delta$-hyperbolic space with bounded growth at some scale is roughly quasi-isometric to a convex subset of $\mathbb{H}^n$. Recently, Wenger \cite{We08} exhibited the largest constant in a quadratic isoperimetric inequality which ensures that a geodesic metric space is Gromov hyperbolic. In \cite{DSU}, Das et al. presented the foundations of the theory of groups and semigroups acting isometrically on Gromov hyperbolic metric spaces. In \cite{LL}, Luo and Kau investigated in detail a simple augmented trees, which was introduced by Kaimanovich \cite{K} on the symbolic spaces of self-similar sets, and the Lipschitz equivalence of such trees. In \cite{HLV07}, Holopainen et al.   studied the asymptotic Dirichlet problem for $p$-harmonic functions in the setting of Gromov hyperbolic metric measure spaces.

\subsection{Uniformity} Uniform domains in Euclidean spaces $\IR^n$ were independently introduced by John in \cite{Jo80} and by Martio and Sarvas in \cite{MS}. The importance of these domains in Euclidean spaces and metric spaces arises from their connections to various results in geometric function theory \cite{BHK,GO,GP,Vai-4}. Since its introduction, researchers have shown an increasing interest in the study of the uniformity from different points of view. For instance, the extendability of functions in Sobolev spaces \cite{Jo81}, geometric characterizations of Gromov hyperbolic domains \cite{BHK,Vai}, Teichm\"uller's displacement problem \cite{MV,VZ}, geometric analysis and potential theory on non-positive curvature spaces \cite{BBS19,KeL,Lo}, and the constructions of planar Loewner carpets in two-dimensional sphere $\mathbb{S}^2$ \cite{CE,EG}.

This paper focuses on the relationship between Gromov hyperbolic spaces and uniform spaces. The connection between these two concepts was observed by Bonk et al. \cite{BHK}  while proving the Gromov hyperbolicity of uniform spaces with respect to the quasihyperbolic metrics. Just as the Poincar\'{e} ball model for $\mathbb{H}^n$, they established a one-to-one correspondence result via replacing the open ball by bounded uniform metric spaces and the space of constant negative curvature by Gromov hyperbolic spaces. Meanwhile, they explained many results in geometric function theory from negative curvature in the sense of Gromov, see also \cite{BHR,Her06,Z}.

\subsection{Uniformity from Gromov hyperbolicity}
It is a natural problem to determine which Gromov hyperbolic domains (that is, an incomplete metric space is $\delta$-hyperbolic with respect to its quasihyperbolic metric) are uniform. This problem has been considered in Euclidean setting such as a bounded domain is uniform if and only if it is $\delta$-hyperbolic and its Euclidean boundary is quasisymmetrically equivalent to the Gromov boundary (\cite[Theorem 1.11]{BHK}), Gromov hyperbolic domains in $\mathbb{S}^2$ are precisely the conformal image of inner uniform slit domains (\cite[Theorem 1.12]{BHK}), a $\delta$-hyperbolic domain in the Riemann sphere with totally disconnected complement is uniform (\cite[Theorem 1.13]{BHK}). There is also a dimension free result in Banach spaces which was proved by V\"{a}is\"{a}l\"{a} in \cite{Vai}.

Along this line, Herron et al. \cite{HSX} generalized their results to metric spaces, and demonstrated that in annular quasiconvex spaces, uniform domains are precisely those Gromov hyperbolic domains whose quasiconformal structure on the Gromov boundary agrees with that on its metric boundary. They also asked whether or not this result holds if the metric space is not annular quasiconvex.

\bques\label{q-1} $($\cite[Question 11.1]{HSX}$)$
Let $(Y,d)$ be a quasiconvex proper metric space and $X\subset Y$ a rectifiably connected open subset with $\partial X\neq \emptyset$. Let $k$ be the quasihyperbolic metric of $X$. Suppose that $(X, k)$ is Gromov hyperbolic, and the natural map exists and is  quasim\"{o}bius. Is $(X,d)$ uniform?
\eques

\begin{remark}\label{r-1}
We recall the definition of natural maps (cf. \cite{BHK,HSX,Vai}). Let $(X,d)$ be an incomplete and rectifiably connected metric space and let $\varrho$ be a metric on $X$ such that $(X,\varrho)$ is Gromov hyperbolic, and let $\partial_\varrho^* X$ and $X_\varrho^*$ be the Gromov boundary and the Gromov closure of $(X,\varrho)$, respectively. For a description of the topology on the Gromov closure, we refer to \cite{BrHa,Vai-0}.

Let $X'$ be the one-point extension $X\cup\{\infty\}$ of $X$ if $X$ is unbounded, and $X'=X$ if $X$ is bounded. Then $\partial X'=\partial X$ if $X$ is bounded, and $\partial X'=\partial X\cup\{\infty\}$ if $X$ is unbounded. If the identity map $\id_X:(X,\varrho)\to (X,d)$ has a continuous extension $\psi:\partial_\varrho^* X\to \partial X'$, then $\psi$ is called a {\it natural} map. Moreover, if such a natural map $\psi$ exists, then we say that $\partial_\varrho^* X$ is naturally equivalent to  $\partial X'$. Furthermore, the mapping $\psi:X_\varrho^*\to \overline{X}'$ with $\psi|_X=\id_X$ is called {\it natural} between the closures, and $X_\varrho^*$ is said to be naturally equivalent to $\overline{X}'$.
\end{remark}

Our  motivation of this paper is to consider Question \ref{q-1}. As the first main result of this article, we construct a counterexample which answers Question \ref{q-1} negatively, see Example \ref{e}.

\begin{thm}\label{thm-1}
There is a compact and geodesic metric space $(Y,d)$ and  a rectifiably connected open subset $X\subset Y$ with $\partial X\neq \emptyset$ such that $(X, k)$ is Gromov hyperbolic, and its Gromov boundary $\partial^*_k X$ equipped with a visual metric is naturally bilipschitz equivalent to the metric boundary $\partial X$. However, $(X,d)$ is not uniform.
\end{thm}

On the other hand, it is interesting to give some affirmative answers to Question $\ref{q-1}$ under certain additional assumptions. We deal with this question in different ways and provide new methods to investigate Gromov hyperbolic characterizations of uniform metric spaces. Our considerations are based on discussions of two hyperbolic type metrics in the sense of Gromov. The first is the quasihyperbolic metric $k$ which was introduced by Gehring and Osgood in \cite{GO}. The second is the hyperbolization metric, denoted by $h$, which was introduced by Ibragimov in \cite{Z1}. Note that Buckley and Herron proved that an incomplete quasiconvex metric space $X$ is uniform if and only if the identity map $(X,k)\to (X,h)$ is rough quasi-isometric, see \cite[Theorem 3.1]{BH}. In view of this result, our approach turns out to be effective to characterize uniform spaces and to study  Question $\ref{q-1}$.

For a bounded incomplete metric space $(X,d)$, we first show that the natural mapping between the Gromov closure of the hyperbolization space $(X,h)$ and the metric completion $\overline{X}$ exists and satisfies a relatively quasisymmetric condition, see Lemma \ref{g-1}. If $X$ is uniform, then a similar result holds for $(X,k)$, see Lemmas \ref{zz3} and \ref{zz4}. So one would like to characterize uniform spaces under the hypothesis  that ``the natural map between the Gromov closure associated with the quasihyperbolic metric and the metric completion of the space is relatively quasisymmetric".

Note that the hyperbolization space $(X,h)$ is not (roughly) geodesic in general. Another useful observation in this paper is that the space $(X,h)$ is {\it BuSc-visual}, see Lemmas \ref{g-3} and \ref{z-2}. This concept was introduced by Buyalo and Schroeder in \cite{BuSc}. We call this property BuSc-visual in order to distinguish another {\it BoS-visual} property on roughly geodesic spaces which was introduced by Bonk and Schramm, see \cite[Section 8]{BS}. The connection between these two  properties are recorded in Lemma \ref{g-5}. It is worth mentioning that the visual properties have lots of applications in recent research such as  in \cite{BHK,HSX,M,ZR,ZLR}.

In view of these observations, we provide a different method to consider the  relationships between uniform metric spaces and Gromov hyperbolic spaces. Note that in \cite{HSX}, the authors pointed out that one cannot obtain a quantitative result even if the answer to Question \ref{q-1} is yes. However, our result will be stated in a quantitative way as follows.

\begin{thm}\label{thm-2} Let $A,C,Q\geq 1$, $K,\delta\geq 0$, $\lambda,\varepsilon\in(0,1)$, and let $\eta:[0,\infty)\to [0,\infty)$ be a homeomorphism. Let $(X,d)$ be a locally compact, rectifiably connected and incomplete metric space. Then $X$ is a bounded $A$-uniform space if and only if the following conditions hold:
\begin{enumerate}
\item[{\rm (1)}] $(X,d)$ is locally $(\lambda,Q)$-quasiconvex;
\item[{\rm (2)}] there is a point $w\in X$ such that $\diam (X)\leq C \dist(w,\partial X)$;
\item[{\rm (3)}] $(X,k)$ is  $\delta$-hyperbolic and $K$-BoS-visual with respect to $w$;
\item[{\rm (4)}] the natural map $\varphi:(X_k^*,k_{w,\varepsilon})\to (\overline{X},d)$ exists and is relatively $\eta$-quasisymmetric.
\end{enumerate}
The constant $A$ and the parameters $C,Q,K,\delta,\lambda,\varepsilon,\eta$ depend only on each other.
\end{thm}

\begin{remark} One observes that Theorem \ref{thm-2} gives a partial answer to Question \ref{q-1} in the affirmative. Under the assumption of that question, the item (1) of Theorem \ref{thm-2} is clear because the quasiconvexity of $Y$ implies that $X\subsetneq Y$ is locally quasiconvex. The item (2) is designed to get a quantitative result, and worth pointing out that there are several examples which show one cannot control the uniformity constant in general, see \cite{HSX,Xie}. It is not difficult to find that the space given in Example \ref{e} is not BoS-visual, which asserts that the item (3) is necessary. The last condition (4) is stronger than that of Question \ref{q-1} where the boundary quasim\"obius equivalence is replaced by a relatively quasisymmetric condition.
\end{remark}

The conformal deformation of a proper geodesic Gromov hyperbolic space introduced by Bonk et al. \cite{BHK} is needed in the proof of Theorem \ref{thm-2}. It turns out that the uniformization theory of Bonk-Heinonen-Koskela allows us to replace Gromov hyperbolic spaces with uniform spaces where the geometry is more easy to understand, and has played a significant role in several research. For example,  Bj\"orn et al. \cite{BBGS17} have regarded the Cantor type sets as boundaries of rooted trees by using this uniformizing theory, and studied the relations between the trace of a first-order Sobolev space and Besov space on the boundary. Applying the conformal density introduced in \cite{BHK}, Carrasco Piaggio \cite{Ca} established a combinatorial description of conformal gauge of a compact space, and showed how to compute its Ahlfors regular conformal dimension. Kemper and Lohkamp in \cite{KeL} and
Lohkamp in \cite{Lo} investigated the relationships between bounded geometry and potential theory for certain elliptic operators on Gromov hyperbolic manifolds through the theory of Bonk-Heinonen-Koskela.

For the unbounded case, with the aid of conformal deformations defined via the Busemann functions on Gromov hyperbolic spaces in \cite{Z}, we obtain the following quantitative result.

\begin{thm}\label{thm-3}   Let $A,Q\geq 1$, $K,\delta\geq 0$, $\lambda,\varepsilon\in(0,1)$, and let $\eta:[0,\infty)\to [0,\infty)$ be a homeomorphism. Let $(X,d)$ be a locally compact, rectifiably connected and incomplete metric space. Then $X$ is an unbounded $A$-uniform space if and only if the following conditions hold:
\begin{enumerate}
\item $(X,d)$ is locally $(\lambda,Q)$-quasiconvex;
\item $(X,k)$ is  $\delta$-hyperbolic and $K$-roughly starlike with respect to $\xi\in\partial_k^*X$;
\item the natural map $\psi:(X_k^*\setminus\{\xi\},k_{b,\varepsilon})\to (\overline{X},d)$ exists with $\psi(\xi)=\infty$ and is relatively $\eta$-quasisymmetric, where the Busemann function $b=b_{\xi,w}\in \mathcal{B}_k(\xi)$ and $w\in X$.
\end{enumerate}
The constant $A$ and the parameters $Q,K,\delta,\lambda,\varepsilon,\eta$ depend only on each other.
\end{thm}

\subsection{Uniformity and Teichm\"{u}ller's  displacement property}\label{sec-t}
As an application of the machinery developed, we characterize uniform metric spaces
whose boundaries are uniformly perfect.

This class of spaces includes all metric spaces $\mathbb{B}_{\rho}^n$ that  arise from a conformal density $\rho$ on the unit balls $\mathbb{B}^n$ (cf.~\cite{BHR,BKR}). It motivates by studying the following problem: {\it For which non-negative measurable function $\rho$ on $\mathbb{B}^n$, can one find a quasiconformal mapping $f$ such that $\rho$ is comparable to a conformal density associated to the Jacobian determinant of $f$?} It should be noted that many results of function theory can be explained through the geometry of $\mathbb{B}^n_\rho$ such as Beurling's radial limit theorem and Gehring-Hayman theorem, see \cite{BKR}. Later, Herron \cite{Her06} generalized the results in \cite{BHR,BKR} to abstract metric measure spaces that is Gromov hyperbolic with uniformly perfect boundaries and bounded geometry. Meanwhile, he proved that for such a space, it admits a uniformizing conformal density if and only if the conformal Assouad dimension of its Gromov boundary is strictly less than its Hausdorff dimension.

As mentioned before, to deal with Question \ref{q-1}, we need additional conditions and consider Teichm\"{u}ller's displacement property for our choice. Let $D\subsetneq \mathbb{R}^n$ be a domain, and let $\mathcal{T}_K(D)$ be the class of quasiconformal homeomorphisms of $D$ onto itself with identity boundary values. Originally, Teichm\"{u}ller \cite{Te} considered the mapping $f\in \mathcal{T}_K(D)$ with $D=\mathbb{R}^2\setminus\{(0,0),(1,0)\}$, and established the following sharp inequality:
$$H_D(x,f(x))\leq \log K$$
for all $x\in D$, where $H_D$ is the classical hyperbolic metric of $D$. For results concerning the same problem in the case of the unit balls, convex domains and uniform domains with uniformly perfect boundaries, we refer to \cite{MV,Martin,Vu1,VZ} and the references therein. In \cite{ZR}, the authors recently studied Teichm\"{u}ller's displacement problem for Gromov hyperbolic domains in $\mathbb{R}^n$ with uniformly perfect Gromov boundaries. We remark that Teichm\"{u}ller type results  are applicable to the theory of homogeneity of domains (cf. \cite{GP,KM11}).

In \cite{Z}, the author investigated Teichm\"{u}ller's displacement problem on Gromov hyperbolic spaces via replacing quasiconformal mappings by rough quasi-isometries. Let $X$ be a proper, geodesic $\delta$-hyperbolic space with $X^*$ its Gromov closure, and let $\mathcal{T}_{M,\mu}(X^*)$ be the class of $(M,\mu)$-roughly quasi-isometric mappings of $X$ onto itself with identity values at the Gromov boundary of $X$. It was shown in \cite{Z} that for all $f\in \mathcal{T}_{M,\mu}(X^*)$, the displacement of each point with its image point under $f$ is bounded above by a uniform constant, provided that the space $X$ is visual and admits uniformly perfect Gromov boundary.

Before stating our results, we need some preparations. Let $(X,d)$ be a rectifiably connected, locally compact and incomplete metric space. In this paper, we consider another class of mappings which is defined as
\beq\nonumber \mathcal{T}_{M,\mu}(\overline{X})
&=& \Big\{{X} \xrightarrow{f} {X}\;:(X,h)\xrightarrow{f}(X,k)\; \mbox{is an $(M,\mu)$-rough quasi-isometry},
\\ \nonumber & & ~~~~~~~\mbox{ the continuous extension of $f$ on the boundary}\; f|_{\partial X}=\id_{\partial X} \Big\}.
\eeq
We say that $X$ satisfies {\it Teichm\"{u}ller's  displacement property} or the {\it displacement property with respect to the quasihyperbolic metric} $k$, if for any constants $M\geq 1,$ $\mu>0$ and for each map $f\in \mathcal{T}_{M,\mu}(\overline{X})$, there exists a function $\Lambda=\Lambda(M,\mu)$ depending only on the data of $f$ and $X$ such that  $k(x,f(x))\leq \Lambda$ for all $x\in X$.

Now, our main result in this direction is as follows.

\begin{thm}\label{thm-4} Let $A,C,Q,\nu\geq 1$, $\delta,\Lambda\geq 0$, $\lambda,\varepsilon\in(0,1)$, and let $\eta:[0,\infty)\to [0,\infty)$ be a homeomorphism. Let $(X,d)$ be a bounded, rectifiably connected, locally compact and incomplete metric space, and let $\partial X$ be $C$-uniformly perfect with $\diam (X)\leq \nu\diam (\partial X)$. Then $X$ is $A$-uniform if and only if the following conditions hold:
\begin{enumerate}
\item  $(X,d)$ is locally $(\lambda,Q)$-quasiconvex;
\item $(X,k)$ is $\delta$-hyperbolic space and the natural map $\varphi:(\partial_k^* X,k_{w,\varepsilon})\to (\partial{X},d)$ exists and is power $\eta$-quasisymmetric;
\item $(X,d)$ satisfies the $\Lambda$-displacement property with respect to the quasihyperbolic metric.
\end{enumerate}
The constant $A$ and the parameters $Q,\Lambda,\delta,\lambda,\varepsilon,\eta$ depend only on each other and $\nu$, $C$.

\end{thm}
For the unbounded case, we obtain the following result which is similar to the above.

\begin{thm}\label{thm-5} Let $A,C,Q\geq 1$, $\delta,\Lambda\geq 0$, $\lambda,\varepsilon\in(0,1)$, and let $\eta:[0,\infty)\to [0,\infty)$ be a homeomorphism. Let $(X,d)$ be an unbounded, rectifiably connected, locally compact and incomplete metric space, and let $\partial X\cup\{\infty\}$ be a $C$-uniformly perfect set. Then $X$ is $A$-uniform if and only if the following conditions hold:
\begin{enumerate}
\item  $(X,d)$ is locally $(\lambda,Q)$-quasiconvex;
\item $(X,k)$ is $\delta$-hyperbolic and the natural map $\psi:(\partial_k^* X\setminus\{\xi\},k_{b,\varepsilon})\to (\partial{X},d)$ exists and is power $\eta$-quasisymmetric, where $\psi(\xi)=\infty$, $b=b_{\xi,w}\in \mathcal{B}_k(\xi)$ and $w\in X$;
\item $(X,d)$ satisfies the $\Lambda$-displacement property with respect to the quasihyperbolic metric.
\end{enumerate}
The constant $A$ and the parameters $Q,\Lambda,\delta,\lambda,\varepsilon,\eta$ depend only on each other and  $C$.
\end{thm}

\begin{remark}
We assume in Theorem \ref{thm-5} that the one-point extended boundary $\partial X\cup\{\infty\}$ of $X$ is $C$-uniformly perfect associated with a sphericalized metric which was introduced by Bonk and Kleiner \cite{BK}. This makes sense because the uniform perfectness is preserved under quasim\"obius maps (cf. \cite{ZR}).
We also remark that Theorems \ref{thm-4} and \ref{thm-5} answer partially to Question \ref{q-1} in the affirmative under the additional Teichm\"uller's displacement property.
\end{remark}

\subsection{Overview of the remainder of the paper}

The rest of this paper is organized as follows. In Section 2, we introduce some notation and necessary terminology concerning metric geometry and mappings on metric spaces.
The focus of Section 3 is Gromov hyperbolic spaces and some auxiliary results about them. In particular, we prove Theorem \ref{thm-1} by presenting a counterexample  in Subsection \ref{b-3}. In Section 4, we introduce the concepts of relative maps and investigate the relationship between rough quasi-isometric maps and relative maps on Gromov hyperbolic spaces. Several techniques and auxiliary lemmas are demonstrated for our need. In Section 5, we investigate the Gromov hyperbolicity of incomplete metric spaces associated with two hyperbolic type metrics, the hyperbolization metric $h$ and the quasihyperbolic metric $k$. Section 6 is devoted to the proofs of Theorems \ref{thm-2} and \ref{thm-3} by using the theory of uniformizing Gromov hyperbolic spaces. The main content of Section 7 is to study the connections between uniform spaces and Teichm\"uller's displacement property. Moreover, the proofs of Theorems \ref{thm-4} and \ref{thm-5} are given.

\section{Preliminaries}

\subsection{Notations} Following \cite{BuSc}, for $s,t,c\in \mathbb{R}$ with $c\geq 0$, it is convenient to write
$$s\doteq t\;\;\mbox{up to an error $\leq c$ or}\;\; s\doteq_c t\;\mbox{instead of}\;\;|s-t|\leq c.$$
Also, we extend our argument about $\doteq$ as follows. For any sequences $\{s_i\},\{t_i\}\subset \mathbb{R}$, we write
$$\{s_i\}_i \doteq \{t_i\}_i\;\;\mbox{up to an error $\leq c$, or}\;\; \{s_i\}_i \doteq_c \{t_i\}_i$$
if $\limsup_{i\to \infty} |s_i-t_i|\leq c$. Similarly, we write
$$\{s_i\}_i\leq \{t_i\}_i+c\;\;\mbox{if}\;\limsup_{i\to \infty} s_i\leq \liminf_{i\to \infty} t_i+c,$$
and
$$\{s_i\}_i\geq \{t_i\}_i+c\;\;\mbox{if}\;\liminf_{i\to \infty} s_i\geq \limsup_{i\to \infty} t_i+c.$$

Moreover, let $A,B>0$ and $C\geq 1$. We use the notation $A\asymp B$ up to multiplicative error $\leq C;$ or $A\asymp_C B$ instead of
$C^{-1}\leq AB^{-1}\leq C.$

For real numbers $s,t$, we write
$$s\wedge t=\min\{s,t\}\;\;\;\;\mbox{and}\;\;\;\; s\vee t=\max\{s,t\}.$$

\subsection{Metric geometry}
Let $(X,d)$ be a metric space. We say that $X$ is {\it incomplete} if its boundary $\partial X= \overline{X}\setminus X\not=\emptyset$, where $\overline{X}$ denotes the metric completion of $X$. For $x\in X$, let $d(x)$ be the distance between $x$ and its boundary $\partial X$.
Let $B(x,r)$  and $\overline{B}(x,r)$ be the open ball and closed ball (of radius $r$ centered at the point $x$) in $(X,d)$, respectively. The space $X$ is called {\it proper} if its closed balls are compact. The diameter of a bounded set $E\subset X$ is denoted by $\diam(E)$.

Let $C\geq 1$, a metric space $X$ is called {\it $C$-uniformly perfect}, if for each $x\in X$ and every $r>0$, $B(x,r)\setminus B(x, r/C)\not=\emptyset$ provided $X\setminus B(x,r)\not=\emptyset$.

A {\it curve} in $X$ is a continuous map $\gamma:\; I\subset \IR\to X$. We also denote the image set $\gamma(I)$ of $\gamma$ by $\gamma$.
If $\gamma$ is an embedding of $I$, then it is called an {\it arc}.  The {\it length} $\ell(\gamma)$ of $\gamma$ with respect to the metric $d$ is defined in an obvious way. Here the parameter interval $I$ is allowed to be open or half-open. We let $\gamma[x,y]$
be the subarc of $\gamma$  between $x$ and $y$ in $\gamma$.  The space $X$ is called {\it rectifiably connected} if every pair of points in $X$ can be connected by a curve $\gamma$ in $X$ with $\ell(\gamma)<\infty$.

A {\it geodesic} $\gamma$ joining $x$ to $y$ in $X$ is a map $\gamma:I=[0,l]\to X$  satisfying $\gamma(0)=x$, $\gamma(l)=y$ and $d(\gamma(t),\gamma(t'))=|t-t'|$ for all $t,t'\in I.$
If $I=[0,\infty)$, then $\gamma$ is called a {\it geodesic ray}. Moreover, if $I=\mathbb{R}$, then $\gamma$ is called a {\it geodesic line}. A metric space $X$ is said to be {\it geodesic} if every pair of points can be joined by a geodesic arc. Let $[x,y]$ denote a geodesic between $x$ and $y$ in $X$.

Let $\varepsilon\geq 0$. Following \cite{BS}, an {\it $\varepsilon$-rough geodesic} is a map $\gamma:I\subset \mathbb{R}\to X$ with
$d(\gamma(t),\gamma(t'))\doteq_{\varepsilon} |t-t'|$, for all $t,t'\in I.$ Note that a rough geodesic may be discontinuous.
The space $X$ is called {\it $\varepsilon$-roughly geodesic} if for every pair of points $x,y\in X$ there exists an $\varepsilon$-rough geodesic connecting $x$ and $y$. Moreover, $X$ is called {\it roughly geodesic} if $X$ is $\varepsilon$-roughly geodesic for some $\varepsilon\geq 0$.

Let $A\geq 1$. A metric space $X$ is called {\it $A$-quasiconvex} if each pair of  points $x,y\in X$ can be joined by an {\it $A$-quasiconvex} curve $\gamma$, that is,
$\ell(\gamma)\leq A d(x,y).$ Let $0<\lambda\leq 1/2$ and $Q\geq 1$. Let $X$ be an incomplete and rectifiably connected metric space. Then
\begin{enumerate}
  \item $X$ is called {\it locally $(\lambda,Q)$-quasiconvex}, if for each $x\in X$ and for any points $y,z\in B(x,\lambda d(x))$, there is a $Q$-quasiconvex curve $\gamma$ in $X$ joining $y$ to $z$.
  \item $X$ is called $A$-{\it uniform} provided  each pair of points $x$, $y$ in $X$ can be joined by an $A$-quasiconvex arc $\gamma$ in $X$ such that
$$\ell(\gamma[x,z])\wedge\ell(\gamma[z,y])\leq A\,d(z) ~\mbox{ for all $z\in \gamma$.}
$$
Also, $\gamma$ is said to be an $A$-{\it uniform} arc. Moreover, we say that $X$ is uniform if it is $A$-uniform for some $A\geq 1$.
  \item The {\it quasihyperbolic metric} $k$ of $X$ is defined by
$$k(x,y)=\inf \int_\gamma \frac{|dz|}{d(z)},$$
where the infimum is taken over all rectifiable curves in $X$ joining $x$ and $y$.
\item The {\it distance ratio metric} is denoted by
$$j(x,y)=\log\left (1+\frac{d(x,y)}{d(x)\wedge d(y)}\right ).$$
\end{enumerate}

Moreover, we record the following useful characterization of uniform spaces.

\begin{thm}\label{zz1}$($\cite[Theorem $3.1$]{BH}$)$ Let $(X,d)$ be an incomplete, rectifiably connected, and locally $(\lambda,Q)$-quasiconvex space. Then $X$  is $A$-uniform if and only if for all $x,y\in X$, we have $k(x,y)\leq c_1j(x,y)+c_2$. The constants $c_1,c_2$ and $A$ depend on each other and $\lambda$ and $Q$.
\end{thm}

\subsection{Mappings on metric spaces}
\bdefe
Let $f:$ $(X,d)\to (Y,d')$ be a mapping (not necessarily continuous) between metric spaces $X$ and $Y$, and let $L\geq 1$ and $\mu\geq 0$ be constants. \begin{enumerate}
\item
If  $L^{-1}d(x,y)-\mu\leq d'(f(x),f(y))\leq Ld(x,y)+\mu$ for all $x,y\in X$,
then $f$ is called an {\it $(L, \mu)$-roughly quasi-isometric mapping} (cf. \cite{BS,HSX}).
If $L=1$, then $f$ is called {\it $\mu$-roughly isometric}.
If in addition
$\sup_{y\in Y}d'(y, f(X))<+\infty,$
then $f$ is called an {\it $(L, \mu)$-rough quasi-isometry}.
\noindent
\item If $f$ is an $(L, 0)$-rough quasi-isometric homeomorphism, then we simply call it as  {\it $L$-bilipschitz}.
\item If  $L d(x,y)-\mu\leq d'(f(x),f(y))\leq Ld(x,y)+\mu$ for all $x,y\in X$,
then $f$ is called {\it $(L, \mu)$-roughly similar}.
\end{enumerate}
\edefe

Two mappings $\varphi,$ $\phi:$ $(X,d)\to (Y,d')$ are said to be {\it roughly equivalent}, denoted by $\varphi\doteq \phi$, if there is a constant  $\mu \geq 0$ such that $d'(\varphi(x),\phi(x))\leq \mu$ for all $x\in X$.
A {\it rough inverse} of a rough quasi-isometry $\varphi:$ $(X,d)\to (Y,d')$ is a rough quasi-isometry $\phi:$ $(Y,d')\to (X,d)$ such that $\phi\circ \varphi\doteq \id_{X}$ and $\varphi\circ \phi\doteq \id_{Y}$. Note that the existence of such a rough inverse can be found in \cite{Vai-0}.

A quadruple in $X$ is an ordered sequence $E=(x,y,z,w)$ of four distinct points
in $X$ and the {\it cross ratio} of $E$ is defined to be the number
$$\tau(E) = \tau(x,y,z,w) = \frac{d(x,z)d(y,w)}{d(x,y)d(z,w)}.$$

Observe that the definition is extended in the well known manner to the case
where one of the points is $\infty$. For example,
$\tau(x,y,z,\infty) =d(x,z)/d(x,y).$

\bdefe A homeomorphism $f:$ $(X,d)\to (Y,d')$ between two metric spaces is said to be
\begin{enumerate}
\item $\eta$-{\it quasisymmetric} if there is a homeomorphism $\eta : [0,\infty) \to [0,\infty)$ such that
$$ d(x,y)\leq td(x,z)\;\; \mbox{implies}\;\;   d'(f(x),f(y)) \leq \eta(t)d'(f(x),f(z))$$
for each $t>0$ and for each  triplet of points $x,$ $y$, $z$  in $X$. If $\eta(t)=C(t^{\alpha}\vee t^{1/\alpha})$ for some $C\geq 1$ and $\alpha\geq 1$, then $f$ is {\it power quasisymmetric}.
\item {\it $\theta$-quasim\"obius} if $\tau(f(E))\leq\theta(\tau(E))$ holds for each quadruple $E\subset X$.
\end{enumerate}\edefe

\section{Gromov hyperbolic spaces}
In this section, we recall necessary terminology concerning Gromov hyperbolic spaces, for which we refer to standard references \cite{BS,BrHa,BuSc,Gr87,LL,Vai-0}.


\subsection{Gromov hyperbolicity}
Let $(X,d)$ be a metric space. Fix a base point $w\in X$. For all $x,y\in X$, we define
$$(x|y)_w=2^{-1}[d(x,w)+d(y,w)-d(x,y)].$$
This number is called the  {\it Gromov product} of $x,y$ with respect to $w$.
We say that $X$ is a {\it Gromov hyperbolic} space, if there is a constant $\delta\geq 0$ such that
$$(x|y)_w\geq (x|z)_w\wedge(z|y)_w-\delta ~\mbox {for all $x,y,z\in X.$}
$$
In other words, the item $((x|y)_w,(x|z)_w,(y|z)_w)$ is a $\delta$-{\it triple}.


\subsection{Hyperbolic metrics and visual metrics based at $w\in X$}
Following \cite{LL}, we define a hyperbolic metric on a Gromov $\delta$-hyperbolic space $(X,d)$ based at some point $w\in X$. For $\varepsilon>0$ with $e^{\varepsilon\delta}\leq 2$, define
$$\rho_{w,\varepsilon}(x,y)=\delta_{x,y}e^{-\varepsilon(x|y)_w}$$
for all $x,y\in X$ with convention $e^{-\infty}=0$, where  $\delta_{x,x}=0$ if $x=y$, and  $\delta_{x,y}=1$ if $x\neq y$. For all $x,y,z\in X$, we have
$\rho_{w,\varepsilon}(x,y)\leq e^{\varepsilon\delta} \big(\rho_{w,\varepsilon}(x,z)\vee \rho_{w,\varepsilon}(z,y)\big).$
That is, $\rho_{w,\varepsilon}$ is a $K$-quasimetric on $X$ with $K=e^{\varepsilon\delta}\leq 2$. We now define
$$\sigma_{w,\varepsilon}(x,y):=\inf \left\{\sum_{i=1}^{n} \rho_{w,\varepsilon} (x_{i-1},x_i):n\geq 1,x=x_0,x_1,\ldots,x_n=y\in X \right\}.$$
According to \cite[Lemma $2.2.5$]{BuSc}, $(X, \sigma_{w,\varepsilon})$ is a metric space with
\be\label{h-1} 2^{-1}\rho_{w,\varepsilon} \leq \sigma_{w,\varepsilon}\leq \rho_{w,\varepsilon}.\ee

Then $\sigma_{w,\varepsilon}$ is called a {\it hyperbolic metric} on $X$ based at $w$ with the parameter $\varepsilon$. The metric completion of $X$ under $\sigma_{w,\varepsilon}$, denoted by $X^*$, is called the {\it Gromov closure} of $X$ and $\partial^* X=X^*\setminus X$ is called the {\it hyperbolic boundary} of $X$. We remark  that the definition of hyperbolic boundary is independent of $w$ and $\varepsilon$ in the sense  that with different choices of base point $w$ and parameter $\varepsilon$, these hyperbolic boundaries are homeomorphic to each other. For more details we refer to \cite{LL}.

An alternative description can be given in terms of {\it the boundary at infinity} or {\it Gromov boundary} $\partial_\infty X$ of $X$, which is  defined to be the set of equivalence classes of sequences converging to infinity; that is, $\{x_i\}\subset X$ {\it converges to infinity}, or say $\{x_i\}$ is a {\it Gromov sequence},  if
$\lim_{i,j\to \infty} (x_i|x_j)_w=\infty.$
So throughout this paper, we take $\partial^* X=\partial_\infty X$ as a set for the simplification of notations. If $(X,d)$ is proper and geodesic, then the hyperbolic boundary is also equivalent to the {\it geodesic boundary} which is defined as the set of equivalence classes of geodesic rays, with two geodesic rays being equivalent if they have finite Hausdorff distance (cf.~\cite{BrHa}).

As usual, the restriction of $\sigma_{w,\varepsilon}$ on the hyperbolic boundary $\partial^*X$ is called a {\it visual metric} based at $w$ with the parameter $\varepsilon$ (cf.~\cite{BS,BuSc}).

\subsection{Hyperbolic metrics and visual metrics based at $\xi\in \partial^*X$}
We begin with the definitions of Gromov products on the hyperbolic completion.

Let $(X,d)$ be a Gromov $\delta$-hyperbolic metric space,  $o\in X$ and  $\xi\in\partial^*X$. For $x\in X$ and $\zeta\in \partial^*X$, the Gromov product $(x|\zeta)_w$ of $x$ and $\zeta$ is defined by
$$(x|\zeta)_w= \inf \left \{ \liminf_{i\rightarrow \infty}(x|z_i)_w:\; \{z_i\}\in \zeta\right\}.
$$
For points $\xi_1,\xi_2\in\partial^*X$, we define their Gromov product by
$$(\xi_1|\xi_2)_w=\inf\left \{\liminf_{i\to\infty} (x_i|y_i)_w: \{x_i\}\in\xi_1 , \{y_i\}\in\xi_2\right \}.
$$

Following \cite{BuSc}, we say that $b:X\to \mathbb{R}$ is a {\it Busemann function} based on $\xi$, denoted by $b\in \mathcal{B}(\xi)$, if
$$b(x)=b_{\xi,o}(x)=b_\xi(x,o)=(\xi|o)_x-(\xi|x)_o ~\mbox{ for all $x\in X$}.$$

We next define the Gromov product of $x,y\in X$ based at the Busemann function $b=b_{\xi,o}\in \mathcal{B}(\xi)$ by
$$(x|y)_b=2^{-1}(b(x)+b(y)-d(x,y)).$$
Similarly, for $x\in X$ and $\zeta\in \partial^*X\setminus\{\xi\}$, the Gromov product $(x|\zeta)_b$ of $x$ and $\zeta$ is defined by
$$(x|\zeta)_b= \inf \left\{ \liminf_{i\rightarrow \infty}(x|z_i)_b:\; \{z_i\}\in \zeta\right\}.
$$
For points $\xi_1,\xi_2\in \partial^*X\setminus\{\xi\}$, we define their Gromov product based at $b$ by
$(\xi_1|\xi_2)_b=\inf\big\{\liminf_{i\to\infty} (x_i|y_i)_b: \{x_i\}\in\xi_1 , \{y_i\}\in\xi_2\}.$
Note that $(\xi_1|\xi_2)_b=-\infty$ if and only if $\xi_1=\xi$ or $\xi_2=\xi$. For more discussions we refer to \cite[Section 3.2.2]{BuSc}.

It follows from \cite[$(3.2)$ and Example $3.2.1$]{BuSc}  that
\be\label{z0} (x|y)_b\doteq_{10\delta} (x|y)_o-(\xi|x)_o-(\xi|y)_o.\ee
Moreover, we observe from \cite[Proporsition $3.2.3$]{BuSc} that  $(x|y)_b,(y|z)_b,(x|z)_b$ form a $22\delta$-triple for every $x,y,z\in X$.

In view of these facts, we introduce a hyperbolic metric on $X$ based at $\xi$ or a Busemann function $b=b_{\xi,w}\in \mathcal{B}(\xi)$ in the same way as the preceding subsection. For $\varepsilon>0$ with $e^{22\varepsilon\delta}\leq 2$, define
$$\rho_{b,\varepsilon}(x,y)=\delta_{x,y}e^{-\varepsilon(x|y)_b}\;\;\;\;\;\;\;\mbox{for}\;x,y\in X.$$
For each $x,y,z\in X$, we have
$\rho_{b,\varepsilon}(x,y)\leq e^{22\varepsilon\delta} \big(\rho_{b,\varepsilon}(x,z)\vee \rho_{b,\varepsilon}(z,y)\big).$
That is, $\rho_{b,\varepsilon}$ is a $K'$-quasimetric on $X$ with $K'=e^{22\varepsilon\delta}\leq 2$. We now define
$$\sigma_{b,\varepsilon}(x,y):=\inf\left\{\sum_{i=1}^{n} \rho_{b,\varepsilon} (x_{i-1},x_i):n\geq 1,x=x_0,x_1,\ldots,x_n=y\in X\right\}.$$
Again by \cite[Lemma $2.2.5$]{BuSc}, $(X, \sigma_{b,\varepsilon})$ is a metric space with
\be\label{h-2} \rho_{b,\varepsilon}/2\leq \sigma_{b,\varepsilon}\leq \rho_{b,\varepsilon}.\ee

We say that $\sigma_{b,\varepsilon}$ is a {\it hyperbolic metric} on $X$ based at $\xi$ or Busemann function $b$ with the parameter $\varepsilon$. The metric completion of $X$ under $\sigma_{b,\varepsilon}$ is denoted by $X_*$ and $\partial_* X=X_*\setminus X$.
Moreover, one easily sees that for all $x_i\neq x_i'\in X$, $\sigma_{b,\varepsilon}(x_i,x_i')\to 0$ if and only if $(x_i|x_i')_b\to \infty$ as $i \to\infty$. According to \cite[Proposition $3.4.1$]{BuSc}, there are canonical inclusions
$$\partial^*X\setminus \{\xi\}\longrightarrow \partial_*X\longrightarrow \partial^*X$$
whose composition coincides with the inclusion $\partial^*X\setminus \{\xi\} \subset \partial^*X$. Thus, the restriction of $\sigma_{b,\varepsilon}$ on the punctured hyperbolic boundary
$\partial^*X\setminus \{\xi\}$ is called a {\it visual metric} based at $\xi$ or the Busemann function $b$ with the parameter $\varepsilon$ (cf. \cite{BuSc}). Note that the visual metric $\sigma_{b,\varepsilon}$ has a natural extension:  $\sigma_{b,\varepsilon}(\xi_1,\xi_2)=\infty$ if and only if $\xi_1=\xi$ or $\xi_2=\xi$, see \cite[Section 3.3]{BuSc}.

Finally, we need some auxiliary results for later use.

\begin{lem}\label{b-1}$($\cite[Lemma $5.11$]{Vai-0}$)$
Let $o,z\in X$, let $X$ be a $\delta$-hyperbolic space, and let $\xi,\xi'\in\partial_\infty X$. Then for any sequences $\{y_i\}\in \xi$, $\{y_i'\}\in \xi'$, we have
\begin{enumerate}
\item  $(z|\xi)_o\leq  \liminf_{i\rightarrow \infty} (z|y_i)_o \leq  \limsup_{i\rightarrow \infty} (z|y_i)_o\leq (z|\xi)_o+\delta;$
\item  $(\xi|\xi')_o\leq  \liminf_{i\rightarrow \infty} (y_i|y_i')_o \leq  \limsup_{i\rightarrow \infty} (y_i|y_i')_o\leq (\xi|\xi')_o+2\delta.$
\end{enumerate}
\end{lem}

\begin{lem}\label{z1}  Let $X$ be a Gromov $\delta$-hyperbolic space and $o,w\in X$, $\xi\in\partial^*X$, and let $\varepsilon, \varepsilon'>0$ with $e^{\delta\varepsilon}\leq 2$ and $e^{22\delta\varepsilon'}\leq2$. Then
\begin{enumerate}
\item  $(X^*,\sigma_{o,\varepsilon})$ is bounded with diameter less than $1$ and $(X_*,\sigma_{b,\varepsilon'})$ is unbounded with $\sigma_{b,\varepsilon'}(\xi,\xi')=\infty$ for any $\xi'\neq \xi$, where $\xi'\in\partial^*X$ and $b=b_{\xi,w}\in \mathcal{B}(\xi)$.
\item  the identity map $(X^*\setminus\{\xi\},\sigma_{o,\varepsilon})\to (X^*\setminus\{\xi\},\sigma_{b,\varepsilon'})$ is $\theta$-quasim\"obius with $\theta(t)=Mt^{\varepsilon'/\varepsilon}$, where $M=M(\delta,\varepsilon, \varepsilon')$.
\end{enumerate}
\end{lem}
\bpf $(1)$ For all $x,y\in X$, by (\ref{h-1}), we have
$\sigma_{o,\varepsilon}(x,y)\leq \rho_{o,\varepsilon}(x,y)\leq e^{-\varepsilon (x|y)_o}\leq 1,$
which guarantees that $\diam(X,\sigma_{o,\varepsilon})\leq 1.$ For any sequences $\{x_n\}\in \xi$ and $\{y_n\}\in \xi'$, we see from $(\ref{z0})$ that for all $n\in \mathbb{N}:=\{1, 2, \ldots\}$,
$$(x_n|y_n)_b\doteq_{10\delta} (x_n|y_n)_w-(\xi|x_n)_w-(\xi|y_n)_w.$$
This implies that $(x_n|y_n)_b\to-\infty$ as $n\to\infty$, because $\xi\neq \xi'$. Thus it follows from Lemma \ref{b-1} that $(\xi|\xi')_b=-\infty$ and hence
$\sigma_{b,\varepsilon'}(\xi,\xi')\geq 1/2e^{-\varepsilon'(\xi|\xi')_b}=\infty,$
the second statement follows.

$(2)$ For any distinct points $x_i\in X^*\setminus\{\xi\}$, $i=1,2,3,4$, we may assume without loss of generality that both of them lie in $X$ by using Lemma \ref{b-1}. Again by $(\ref{z0})$, an elementary computation gives that

\vspace{8pt}
\noindent
$(x_1|x_2)_o+(x_3|x_4)_o-(x_1|x_3)_o-(x_2|x_4)_o$
\begin{align}\nonumber
& = \;\;\;\; (x_1|x_2)_w+(x_3|x_4)_w-(x_1|x_3)_w-(x_2|x_4)_w
\\ \nonumber
&\doteq_{40\delta} (x_1|x_2)_b+(x_3|x_4)_b-(x_1|x_3)_b-(x_2|x_4)_b.
\end{align}
Now by (\ref{h-1}) and (\ref{h-2}), we compute
\beq\nonumber
\frac{\sigma_{b,\varepsilon'}(x_1,x_3)\sigma_{b,\varepsilon'}(x_2,x_4)}{\sigma_{b,\varepsilon'}(x_1,x_2)\sigma_{b,\varepsilon'}(x_3,x_4)}&\leq& 4e^{-\varepsilon'[(x_1|x_2)_b+(x_3|x_4)_b-(x_1|x_3)_b-(x_2|x_4)_b]}
\\ \nonumber &\leq& 4e^{-40\varepsilon'\delta}e^{-\varepsilon'[(x_1|x_2)_o+(x_3|x_4)_o-(x_1|x_3)_o-(x_2|x_4)_o]}
\\ \nonumber &\leq& 4e^{-40\varepsilon'\delta}4^{\varepsilon'/\varepsilon}
\Big[\frac{\sigma_{o,\varepsilon}(x_1,x_3)\sigma_{o,\varepsilon}(x_2,x_4)}{\sigma_{o,\varepsilon}(x_1,x_2)
\sigma_{o,\varepsilon}(x_3,x_4)}\Big]^{\varepsilon'/\varepsilon},
\eeq
as desired. This completes the proof.
\epf

\subsection{Visual properties of Gromov hyperbolic spaces} This part focuses on the relations between two visual properties of $\delta$-hyperbolic spaces which were independently introduced by Bonk and Schramm \cite{BS} and by Buyalo and Schroeder \cite{BuSc}. We start with the following definitions.

%

\bdefe Let $(X,d)$ be a Gromov hyperbolic space with $w\in X$, and $K\geq 0$. Then $X$ is said to be
\begin{enumerate}
\item    $K$-{\it BoS-visual} with respect to $w$ if every $x\in X$ lies in a $K$-roughly geodesic ray emanating from $w$.

\item  $K$-{\it BuSc-visual} with respect to $w$ if for every $x\in X$, there is a $\xi\in\partial^*X$ which satisfies
$d(w,x)\leq (x|\xi)_{w}+K.$
\end{enumerate}

Moreover, we say that $X$ is {\it BoS-visual} (resp {\it BuSc-visual}) if it is $K$-{\it BoS-visual} (resp. $K$-{\it BuSc-visual})  with respect to $w$ for some $K\geq 0$ and $w\in X$.
\edefe

We remark that the $BoS$-visual property is equivalent to the rough starlikeness of a $\delta$-hyperbolic space $X$ provided it is proper and geodesic. Following \cite{BHK}, $X$ is called $K$-{\it roughly starlike} with respect to $w\in X$ if for each $x\in X$, there is some point $\xi\in\partial^*X$ and a geodesic ray $\alpha$ connecting $w$ to $\xi$ such that $\dist(x,\alpha)\leq K$.

Next, we recall the definition of the rough starlikeness with respect to a boundary point at infinity.
\bdefe
Let $X$ be a proper, geodesic $\delta$-hyperbolic space, and $\xi\in\partial^* X$. We say that $X$ is {\it $K$-roughly starlike} with respect to $\xi$ if for each $x\in X$ there is a point $\zeta\in\partial^*  X$ and a geodesic line $\gamma=[\xi,\zeta]$ joining $\xi$ and $\zeta$ with $\dist(x,\gamma)\leq K$.
\edefe

Note that every $\delta$-hyperbolic domain in $\mathbb{R}^n$ is roughly starlike with respect to each point of the boundary at infinity, see \cite[Theorem $3.22$]{Vai}. In metric spaces setting, we have the following:

\begin{lem}\label{z2}$($\cite[Lemma 3.14]{Z}$)$  Let $(X,d)$ be an unbounded, locally compact $A$-uniform space. Then there is a natural identification between $\partial^*X$ and $\partial X\cup\{\infty\}$ and
$(X,k)$ is $K$-roughly starlike with respect to each point of $\partial^* X$ with $K=K(A)$.
\end{lem}

For a roughly geodesic hyperbolic space $X$, we show that the above two visual properties are equivalent. Before the proof, we need some preparation.

\begin{lem}\label{g-4}  Let $X$ be a $\delta$-hyperbolic space and let $\gamma\subset X$ be an $\varepsilon$-rough geodesic with endpoints $x$ and $y$. Then for all $p\in X$, we have
\be\label{h-00} (x|y)_p-\frac{3\varepsilon}{2}\leq \dist(p,\gamma)\leq (x|y)_p+2\delta+\frac{11\varepsilon}{2}. \ee
\end{lem}

\bpf Let $\gamma: [0,s]\to X$ with $\gamma(0)=x$ and $\gamma(s)=y$. For all $t\in [0,s]$, we have
\be\label{gg-0} d(x,y)=d(\gamma(0),\gamma(s))\geq s-\varepsilon\geq d(\gamma(0),\gamma(t))+d(\gamma(t),\gamma(s))-3\varepsilon.\ee
This guarantees that
$\dist(p,\gamma(t))-(x|y)_p\geq -3\varepsilon/2,$
which gives the lower bound estimate of (\ref{h-00}).

For the other direction, we first check that for each $t\in [0,s]$,
\be\label{h-01} (x|p)_{\gamma(t)}\wedge (p|y)_{\gamma(t)}\leq \delta+\frac{3\varepsilon}{2}. \ee
Since $X$ is $\delta$-hyperbolic, we obtain by (\ref{gg-0})
\begin{eqnarray*} (x|p)_{\gamma(t)}\wedge (p|y)_{\gamma(t)} &\leq& (x|y)_{\gamma(t)}+\delta
\\ \nonumber&=& \frac{1}{2}\big(d(\gamma(0),\gamma(t))+d(\gamma(t),\gamma(s))-d(\gamma(0),\gamma(s))\big)+\delta
\\ \nonumber&\leq& \delta+\frac{3\varepsilon}{2},
\end{eqnarray*}
and thus, (\ref{h-01}) is true. In the following, we consider two cases.

For the first case where $(p|y)_x\leq \varepsilon$, we have
$$\varepsilon\geq (p|y)_x=d(p,x)-(x|y)_p\geq \dist(p,\gamma)-(x|y)_p,$$
as required.

For the remaining case that $(p|y)_x> \varepsilon$, we note that
$$t_0:=(p|y)_x-\varepsilon \leq d(x,y)-\varepsilon\leq s$$
and $\gamma(t_0)=q$. We compute
\begin{eqnarray*} |2(p|y)_q-2(p|x)_q| &=& |d(\gamma(t_0),\gamma(s))-d(\gamma(t_0),\gamma(0))+d(p,x)-d(p,y)|
\\ \nonumber&\leq& |s-2t_0+d(p,x)-d(p,y)|+2\varepsilon
\\ \nonumber&\leq& |d(x,y)-2(p|y)_x+d(p,x)-d(p,y)|+5\varepsilon=5\varepsilon.
\end{eqnarray*}
This, together with $(\ref{h-01})$, shows that
\be\label{h-02} (p|y)_q+(p|x)_q \leq 2 ((x|p)_q\wedge (p|y)_q)+\frac{5\varepsilon}{2}\leq 2\delta+\frac{11\varepsilon}{2}.\ee
On the other hand, we observe that
\begin{eqnarray*}(p|y)_q+(p|x)_q &=& d(p,q)+\frac{1}{2}\big(d(y,q)+d(x,q)-d(p,y)-d(p,x)\big)
\\ \nonumber&\geq& \dist(p,\gamma)-(x|y)_p.
\end{eqnarray*}
Thus by (\ref{h-02}), we obtain the upper bound estimate of (\ref{h-00}). The proof of Lemma \ref{g-4} is complete.
\epf

\begin{lem}\label{g-5}  Let $X$ be a $\delta$-hyperbolic, $\varepsilon$-roughly geodesic space. Then $X$ is BuSc-visual if and only if $X$ is BoS-visual.
\end{lem}

\bpf For the necessity, we assume that $X$ is BuSc-visual. It follows that there is a constant $K_1\geq 0$ and a base point $o\in X$ such that for each $x\in X$, there is a point $\xi\in \partial^* X$ satisfying
$d(o,x)\leq (x|\xi)_o+K_1.$
By \cite[Proposition 5.2]{BS}, we see that there is an $\varepsilon_0$-roughly geodesic ray $\gamma:[0,\infty)\to X$ emanating from $\gamma(0)=o$ to $\gamma(\infty)=\xi$ with $\varepsilon_0=\varepsilon_0(\varepsilon,\delta).$ Moreover, for all $t>0$, it follows that $\gamma|_{[0,t]}:=\gamma_t\subset X$ is an $\varepsilon_0$-rough geodesic segment.

By Lemma \ref{b-1}, without loss of generality we may assume that for $t$ sufficiently large,
$d(o,x)\leq (x|\xi)_o+K_1 \leq (x|\gamma(t))_o+K_1,$
which implies
$(o|\gamma(t))_x\leq K_1.$
Using  (\ref{h-00}), we obtain
$$\dist(x,\gamma_t)\leq  (o|\gamma(t))_x+2\delta+6\varepsilon\leq K_1+2\delta+6\varepsilon.$$
We may assume that $d(x,\gamma(t_0))\leq  K_1+2\delta+6\varepsilon$ for some $t_0\in [0,t]$. Let
$$\alpha(t)=\begin{cases}
\displaystyle \;\;\;x \;\;\;\;\;\;\;\;
\mbox{if}\;\; t=t_0,\\
\displaystyle \;\,\;\gamma(t)
\;\;\;\;\mbox{otherwise}.
\end{cases}
$$
It is not hard to see that $\alpha$ is a $K_2$-roughly geodesic with $K_2=K_1+2\delta+6\varepsilon+\varepsilon_0$ and thus, the necessity follows.

Next, we show the sufficiency. Fix $x\in X$. Thus there is a base point $o\in X$ and a nonnegative constant $K$ such that $x$ lies in a $K$-roughly geodesic ray $\beta:[0,\infty)\to X$  emanating from $o=\beta(0)$ with $\beta(\infty)=\xi$. Set $x=\beta(t_0)$. For all $t\geq t_0$, it follows that
$$ 2d(x,o)-2(x|\beta(t))_o = d(\beta(t_0),\beta(0))-d(\beta(t),\beta(0))+d(\beta(t_0),\beta(t))\leq 3K.$$
By letting $t\to \infty$, Lemma \ref{b-1} ensures that
$$d(x,o)\leq (x|\xi)_o +2K+2\delta.$$
Therefore, $X$ is $K_0$-BuSc-visual with $K_0=2K+2\delta$.
The lemma follows.
\epf

\subsection{Proof of Theorem \ref{thm-1}}\label{b-3} The aim of this subsection is to present a counterexample to show Theorem \ref{thm-1}.

\begin{exam}\label{e} Let
$$Y=\big( [0,1]\times\{0\} \big)\cup \bigcup_{n=2}^{\infty} \left(\Big\{\frac{1}{n}\Big\}\times \bigg[0,\frac{1}{\sqrt{n}}\bigg] \right)\subset \mathbb{R}^2,$$
and let $d$ be the length metric on $Y$ induced by the Euclidean metric of $\mathbb{R}^2$.
Write $o=(0,0), a=(1,0), x_n=(1/n,1/\sqrt{n})$ and $y_n=(1/n,0)$ for $n=2,3,\ldots$. Set $X=Y\setminus\{o,a\}$. Then we have the following:
\begin{enumerate}
\item $(Y,d)$ is a compact and geodesic metric space;
\item $(X,k)$ is Gromov $\delta$-hyperbolic with $\delta=0$;
\item the hyperbolic boundary $(\partial_k^*X,k_{y_2,\varepsilon})$ is naturally bilipschitz equivalent to $(\partial X,d)$;
\item $(X,d)$ is not uniform.
\end{enumerate}
\end{exam}

\begin{figure}[htbp]
\begin{center}
\includegraphics{example-3_8.1}
\end{center}
\caption{Example \ref{e}} \label{1}
\end{figure}

\bpf $(1)$ It follows from the construction that the comb space $(Y,d)$ is a compact metric tree, and thus
$(Y, d)$ is a unique geodesic metric space; note that a metric space is called unique geodesic, provided that the geodesic connecting any pair of points in the space exists and is unique.

$(2)$ According to \cite[Proposition $2.8$]{BHK}, $(X,k)$ is a proper geodesic metric space, where $k$ is the quasihyperbolic metric of $X$ with respect to $d$. For every $x,y\in X$, the unique $d$-geodesic $\alpha$ joining $x$ to $y$ in $X$ is a quasihyperbolic geodesic connecting these two points. Moreover, it is easy to see that $(X,k)$ is a unique geodesic metric space and so is a metric tree. Hence $(X,k)$ is Gromov $\delta$-hyperbolic with $\delta=0$.

$(3)$ Firstly, we observe that the intervals $[y_2,a)$ and $[y_2,o)$ are quasihyperbolic geodesic rays in $(X, k)$, and any other quasihyperbolic geodesic rays in $(X, k)$ should be equivalent to one of them. Moreover, $[y_2,a)$ ends at a unique point $a\in \partial X$ and similarly, $[y_2,o)$ ends at a unique point $o\in \partial X$. Since $\partial X=\{o,a\}$, it follows that there is a bijection
$\varphi: \partial_k^* X\to\partial X$
with $\varphi(\overline{a})=a$ and $\varphi(\overline{o})=o$.

Next we prove that
\bcl\label{z-0.3} $\varphi$ is a natural map. \ecl

It suffices to show that every Gromov sequence $\{a_n\}$ in $(X,k)$ $d$-converges to a point on $\partial X$. Put $a_n=(t_n,s_n)$. We have two cases at hands.

If $\{a_n\}\in \overline{a}$, then $s_n=0$ and $t_n>1/2$ for $n$ sufficiently large. Moreover, we have
$$\log\frac{1}{2(1-t_n)}=k(y_2,a_n)=(a_n|a_n)_{y_2}\to \infty \;\;\;\mbox{as}\;n\to \infty,$$
which implies that $t_n\to 1$. Hence $a_n\to a=\varphi(\overline{a})$.

If $\{a_n\}\in \overline{o}$, then with no loss of generality in assuming that $t_n<1/2$. For each $n$, there is an integer $m\geq 2$ (depending on $n$) such that
$\frac{1}{1+m}\leq t_n< \frac{1}{m}.$
An elementary computation shows that
$$\infty\leftarrow (a_n|a_n)_{y_2}=k(y_2,a_n)\leq k(x_{m+1},y_2)\leq 2\log(1+m),$$
and therefore, $m\to \infty$ and $t_n\to 0$ as $n\to \infty.$ Since
$$d(a_n,o)\leq d(o,y_m)+d(x_{1+m},y_{1+m})= \frac{1}{m}+\frac{1}{\sqrt{1+m}},$$
we obtain $a_n\to o=\varphi(\overline{o})$. Hence Claim \ref{z-0.3} is true.

To prove $(3)$, we only need to show that $\varphi$ is bilipschitz. On the one hand, we have $\partial X=\{a,o\}$ and $d(o,a)=1$.
On the other hand, since $(X,k)$ is Gromov $0$-hyperbolic, it is boundary continuous in the sense of \cite[Section $3.4.2$]{BuSc}. For any sequence $\{b_n\}\in \overline{o}\in\partial_k^*X$ and $\{a_n\}\in \overline{a}\in \partial_k^* X$, we have
$(a_n|b_n)_{y_2}\to (\overline{a}|\overline{o})_{y_2}$ as $n\to\infty.$
Choose a Gromov sequence $\{b_n\}\subset[0,1]\times\{0\}$ with $\{b_n\}\in \overline{o}$. Similarly, choose another Gromov sequence $\{a_n\}\subset[0,1]\times\{0\}$ with $\{a_n\}\in \overline{a}$. It follows that for any $n\in  \mathbb{N}_0:=\mathbb{N}\cup\{0\}$, $(a_n|b_n)_{y_2}=0.$
This is because the union of the quasihyperbolic geodesics $[a_n,y_2]_k$ and $[y_2,b_n]_k$ is a quasihyperbolic geodesic connecting $a_n$ and $b_n$. Moreover, by Lemma \ref{b-1}, we have
$0\leq (\overline{a}|\overline{o})_{y_2}\leq \lim_{n\to\infty}(a_n|b_n)_{y_2}=0$
and thus,
$$k_{y_2,\varepsilon}(\overline{o},\overline{a}) \asymp_C e^{-\varepsilon (\overline{a}|\overline{o})_{y_2}}=1=d(o,a).$$
This shows that the natural map $(\partial_k^*X,k_{y_2,\varepsilon})\to (\partial X,d)$ exists and is bilipschitz.

$(4)$  Suppose to the contrary there is a positive number $A$ such that $(X,d)$ is $A$-uniform. By \cite[Theorem $2.10$]{BHK}, for any $n=3,4,...$, we see that the quasihyperbolic geodesic $[x_n,x_2]_k$, which is also a $d$-geodesic, is $B$-uniform for some constant $B=B(A)>0$. Therefore, we have
$$1/\sqrt{n}=d(x_n,y_n)= d(x_n,y_n)\wedge d(y_n,x_2)\leq Bd(y_n)=B/n,$$
which implies $B\geq \sqrt{n}\to \infty$ as $n\to\infty$. Contradiction.
\epf

\section{Relative maps on Gromov hyperbolic spaces}
In this section, we introduce the concepts of relative maps and investigate the relationships between roughly quasi-isometric maps and relative maps on Gromov hyperbolic spaces. In particular, several auxiliary lemmas are required for our discussion.

\subsection{Cross-difference and strong PQ-isometric maps} In this part, we recall the definitions of cross-difference and strong PQ-isometric maps. Note that these concepts were introduced by Buyalo and Schroeder \cite[Chapter 4]{BuSc} in the study of asymptotic behavior of rough quasi-isometries on Gromov hyperbolic spaces.

Let $(X,d)$ be a metric space and $(x,y,z,u)\subseteq X$. The classical {\it cross-difference} is given by
$$\langle x,y,z,u\rangle=\frac{1}{2}\big(d(x,z)+d(y,u)-d(x,y)-d(z,u)\big).$$

We say that a map $f:$ $(X,d)\to (Y,d')$ between metric spaces is {\it strongly PQ-isometric} if there are constants $c_1\geq 1$ and $c_2\geq 0$ such that for all quadruples $(x,y,z,u)\subseteq X$ with $\langle x,y,z,u\rangle\geq 0$, we have
$$c_1^{-1}\langle x,y,z,u\rangle-c_2 \leq \langle f(x),f(y),f(z),f(u)\rangle \leq c_1\langle x,y,z,u\rangle +c_2.$$

A direct computation shows that $\langle x,y,z,u\rangle=-\langle x,z,y,u\rangle$. Thus it is possible to write the condition to be a strongly PQ-isometric map as
$$-\theta\big(-\langle x,y,z,u\rangle\big)\leq \langle f(x),f(y),f(z),f(u)\rangle\leq \theta\big(\langle x,y,z,u\rangle\big),$$
where now  $(x,y,z,u)\subseteq X$ is an arbitrary quadruple and $\theta:\mathbb{R}\to \mathbb{R}$ is the control function
$\theta(t)=(ct)\vee(t/c)+c'$.

\begin{lem}\label{z3} Let $f:(X,d)\to (Y,d')$ be a map between metric spaces.
\begin{enumerate}
\item If $f$ is roughly isometric, then $f$ is strongly $(M,c')$-PQ-isometric;
\item If $f$ is $(M,\mu)$-roughly similar, then $f$ is strongly $(c,c')$-PQ-isometric with $c=M$;
\item $($\cite[Theorem $4.4.1$]{BuSc}$)$ If $X$ and $Y$ are geodesic, $\delta$-hyperbolic spaces and $f$ is $(M,\mu)$-roughly quasi-isometric, then $f$ is strongly $(c,c')$-PQ-isometric with $c$ and $c'$ depending only on $M,\mu$ and $\delta$.
\end{enumerate}
\end{lem}
\subsection{Relative maps}
Let $M\geq 1$, $\alpha>0$ and $\theta:[0,\infty)\to[0,\infty)$ be a homeomorphism.
Let $(X,d)$ and $(Y,d')$ be incomplete metric spaces and  $f:(\overline{X},d)\to (\overline{Y},d')$ be a map such that $f:\partial X\to \partial Y$ is a nonconstant surjective. Then we say that $f$ is
\begin{enumerate}
\item {\it relatively $M$-bilipschitz} if
$$M^{-1}d(x,y)\leq d'(f(x),f(y))\leq Md(x,y)\ \mbox{ for $(x,y)\in [\overline{X}\times \overline{X}] \setminus  [X\times X]$}.
$$

\item {\it relatively $(M,\alpha)$-snowflake} if
$$M^{-1}d(x,y)^\alpha \leq d'(f(x),f(y))\leq Md(x,y)^\alpha \ \mbox{ for $(x,y)\in [\overline{X}\times \overline{X}] \setminus  [X\times X]$}.$$

\item {\it relatively $\theta$-quasisymmetric} if for all triple $(x,y,z)\subset \overline{X}$ with $x\in\partial X$ or $\{y,z\}\subset \partial X$, we have
\be\label{z4} d'(f(x),f(y))\leq \theta(t)d'(f(x),f(z)) ~\mbox{ whenever}\;d(x,y)=td(x,z).
\ee
If $\theta(t)=C(t^{\beta}\vee t^{1/\beta})$ for some $C\geq 1$ and $\beta\geq 1$, then $f$ is called {\it relatively power quasisymmetric}.

\item {\it relatively $\theta$-quasim\"{o}bius} if  for all quadruple $E=(x,y,z,u)\subset \overline{X}$ with $\{x,z
    \}\subset \partial X$ or $\{y,u\}\subset \partial X$, we have
    $\tau(f(E))\leq\theta(\tau(E)).$
\end{enumerate}

\begin{remark} $(1)$ It follows from the definitions that relatively bilipschitz maps are relatively snowflake and relatively snowflake maps are relatively quasisymmetric. The composition of relatively bilipschitz maps are also relatively bilipschitz. Similar assertions hold true for relatively snowflake maps and relatively quasisymmetric maps.

$(2)$  Carefully checking the above definitions, it is not difficult to see that relatively bilipschitz maps and relatively snowflake maps  are both homeomorphisms on the boundary. According to \cite[Theorem $2.21$]{TV}, the restriction of relatively quasisymmetric maps on the boundary are quasisymmetric homeomorphisms. Moreover, to verify the relative quasisymmetry, a similar argument as in \cite[Lemma $6.17$]{Vai-4} ensures that we only need to show the inequality $(\ref{z4})$ for each triple $(x,y,z)\subset \overline{X}$ with $x\in\partial X$.
\end{remark}

\begin{lem}\label{z-1} Suppose that $X$ is unbounded and that $f:(\overline{X},d)\to (\overline{Y},d')$ is relatively quasim\"obius. Then $f$ is relatively quasisymmetric if and only if $f(x)\rightarrow\infty$ as $x\rightarrow \infty$.
\end{lem}
\bpf If $f$ is relatively quasisymmetric, then $f$ is quasisymmetric on $\partial X$. Thus by \cite[Theorem $3.10$]{Vai-0}, we find that $f(x)\rightarrow\infty$ as $x\rightarrow \infty$. On the other hand, since $\tau(x,y,z,\infty)=d(z,x)/d(y,x)$, the assertion follows.
\epf

In the remaining part of this section, we assume that $(X,d)$ is a $\delta$-hyperbolic space with $w\in X$ and $\xi\in\partial^*X$. Let $\sigma_{w,\varepsilon}$ be a hyperbolic metric on $X^*$ based at $w$ with the parameter $\varepsilon=\varepsilon(\delta)$, and $\sigma_{b,\varepsilon}$ be a hyperbolic metric on $X_*=X^*\setminus\{\xi\}$ based at the Busemann function $b\in \mathcal{B}(\xi)$ with the parameter $\varepsilon$. Similarly, $(Y,d')$ is a $\delta$-hyperbolic space with $w'\in Y$ and $\xi'\in\partial^*Y$. Let $\sigma_{w',\varepsilon'}$ be a hyperbolic metric on $Y^*$ based at $w'$ with the parameter $\varepsilon'=\varepsilon'(\delta)$, and $\sigma_{b',\varepsilon'}$ be a hyperbolic metric on $Y_*=Y^*\setminus\{\xi'\}$ based at the Busemann function $b'\in \mathcal{B}(\xi')$ with the parameter $\varepsilon'$.

\begin{lem}\label{z5}
Let $f:X\to Y$ be an $(M,\mu)$-roughly similar map with $f(w)=w'$. If $f^*|_{\partial^*X}$, the restriction of the induced map $f^*:(X^*,\sigma_{w,\varepsilon}) \to (Y^*,\sigma_{w',\varepsilon'})$ of $f$ on $\partial^*X$, is surjective, then $f^*$ and $f_*:(X^*\setminus\{\xi\},\sigma_{b,\varepsilon}) \to (Y^*\setminus\{\xi'\},\sigma_{b',\varepsilon'})$ are both relatively snowflake, where $\xi'=f^*(\xi)$.
\end{lem}
\bpf Fix a pair of points $(x,y)\in (X^*\times X^*)\setminus (X\times X)$ with $x\neq y$. For any sequences $\{x_n\}\in x$ and $\{y_n\}\in y$, we may take $x_n=x$ or $y_n=y$ for all $n\in \mathbb{N}$ if $x$ or $y$ belongs to $X$. With not much is lost in assuming that $x_n\neq y_n$ for all $n$ sufficiently large.

We first show that $f^*$ is relatively snowflake. Note that by (\ref{h-1}), we have
$$\sigma_{w,\varepsilon}(x_n,y_n)\asymp_{C(\delta)}e^{-\varepsilon(x_n|y_n)_w}\;\;\;\;\mbox{and}\;\;\;\;  \sigma_{w',\varepsilon'}(f(x_n),f(y_n))\asymp_{C(\delta)} e^{-\varepsilon'(f(x_n)|f(y_n))_{w'}}. $$
Moreover, by Lemma \ref{z3}, it follows that
$$(f(x_n)|f(y_n))_{w'}=\langle f(x_n),f(y_n),f(w),f(w) \rangle \doteq_{2\mu} M \langle x_n,y_n,w,w \rangle = M (x_n|y_n)_w.$$
Thus we obtain
$$\sigma_{w',\varepsilon'}(f(x_n),f(y_n))\asymp_{C(\delta,\mu,\varepsilon,\varepsilon')} \sigma_{w,\varepsilon}(x_n,y_n)^{M\varepsilon'/\varepsilon},$$
which guarantees that
$$\sigma_{w',\varepsilon'}(f^*(x),f^*(y))\asymp_{C(\delta,\mu,\varepsilon,\varepsilon')} \sigma_{w,\varepsilon}(x,y)^{M\varepsilon'/\varepsilon},$$
as desired.

Next, we prove that $f_*$ is relatively snowflake. For any sequence $\{u_n\}\in \xi$, it follows from \cite[Proposition 6.3]{BS} that $\{f(u_n)\}\in \xi'=f^*(\xi)$. By (\ref{h-2}), we get
$$\sigma_{b,\varepsilon}(x_n,y_n)\asymp_{C(\delta)}e^{-\varepsilon(x_n|y_n)_b}\;\;\mbox{and}\;\;  \sigma_{b',\varepsilon'}(f(x_n),f(y_n))\asymp_{C(\delta)} e^{-\varepsilon'(f(x_n)|f(y_n))_{b'}}. $$
Furthermore, we compute

\vspace{8pt}
\noindent
$(f(x_n)|f(y_n))_{b'}$
\begin{align}\nonumber
 &\doteq_{C(\delta)}(f(x_n)|f(y_n))_{w'}-(\xi'|f(x_n))_{w'}-(\xi'|f(y_n))_{w'}\;\;\;\;\;\;\;\;\;\;\;\;\;\;\;\;\;\;\;\;\; \;\; (\mbox{by}\;(\ref{z0}))
\\ \nonumber &\doteq_{C(\delta)} \{(f(x_n)|f(y_n))_{w'}-(f(u_m)|f(x_n))_{w'}-(f(u_m)|f(y_n))_{w'}\}_m \;(\mbox{by Lemma \ref{b-1}})
\\ \nonumber &\doteq_{C(\delta,\mu)}M\{ (x_n|y_n)_w-(u_m|x_n)_w-(u_m|y_n)_w\}_m\;\;\;\;\;\;\;\;\;\;\;\;\;\; \;\;\;\;\;\;\;  (\mbox{by}\;\mbox{Lemma \ref{z3}})
\\ \nonumber &\doteq_{C(\delta,\mu,M)} M (x_n|y_n)_b.\;\;\;\;\;\;\;\;\;\;\;\;\;\;\;\; \;\;\;\;\;\;\;\;\;\;\;\;\;\;\; \;\;  \;\;\;\; \;\;\;\; \;\;\;\; \;\;\;\;\;\;\;\;\;\;\; \;\; \;\; \;\;\;\;\;\;\;\;\;\;\;\;\;(\mbox{by}\;(\ref{z0}))
\end{align}
This implies that
$$\sigma_{b',\varepsilon'}(f(x_n),f(y_n))\asymp_{C(\delta,M,\mu,\varepsilon,\varepsilon')} \sigma_{b,\varepsilon}(x_n,y_n)^{M\varepsilon'/\varepsilon},$$
and therefore,
$$\sigma_{b',\varepsilon'}(f_*(x),f_*(y))\asymp_{C(\delta,\mu,M,\varepsilon,\varepsilon')} \sigma_{b,\varepsilon}(x,y)^{M\varepsilon'/\varepsilon},$$
as required. The proof is complete.
\epf

\begin{remark}
The existence of an extension of $f$ on $X^*$ follows from Lemma \ref{z3}$(2)$ and \cite[Corollary 4.3.3]{BuSc}. If we set $M=1$ in Lemma \ref{z5}, that is, if $f$ is roughly isometric, then the induced maps $f^*$ and $f_*$ are both relatively bilipschitz whenever $\varepsilon=\varepsilon'$.
\end{remark}

\begin{lem}\label{zz} Let $X$ and $Y$ be geodesic spaces, and $f:X\to Y$ be an $(M,\mu)$-roughly quasi-isometric map with $f(w)=w'$. If $f^*|_{\partial^*X}$, the restriction of the induced map $f^*:(X^*,\sigma_{w,\varepsilon}) \to (Y^*,\sigma_{w',\varepsilon'})$ on $\partial^*X$, is surjective, then $f^*$ and $f_*:(X^*\setminus\{\xi\},\sigma_{b,\varepsilon}) \to (Y^*\setminus\{\xi'\},\sigma_{b',\varepsilon'})$ are both relatively quasisymmetric, where $\xi'=f^*(\xi)$.
\end{lem}

\bpf We first show that $f^*$ is relatively quasisymmetric. Because $f^*|_{\partial^*X}$ is surjective, it suffices to verify that there is a homeomorphism $\theta:[0,\infty)\to[0,\infty)$ such that for all distinct points  $\{x,y,z\}\subset X^*$ with $x\in\partial^* X$,
$$T=\frac{\sigma_{w',\varepsilon'}(f^*(x),f^*(y))}{\sigma_{w',\varepsilon'}(f^*(x),f^*(z))}\leq \theta(t)\;\;\;\;\;\;\;\;\mbox{whenever}\;\sigma_{w,\varepsilon}(x,y)=t\sigma_{w,\varepsilon}(x,z).$$

Let $\{x_n\}\in x$ be a Gromov sequence. We assume that $y,z\in \partial^* X$ and take two Gromov sequences $\{y_n\}\in y$ and $\{z_n\}\in z$. (This is because if $y$ or $z$ in $X$, then we take $y_n=y$ or $z_n=z$ for all $n\in \mathbb{N}$. The argument below follows in a natural way.)
Without loss of generality we assume that any two points of $\{x_n, y_n,z_n\}$ are not the same for all $n$ sufficiently large. Using \cite[Proposition 6.3]{BS}, one finds that $\{f(x_n)\}\in f^*(x)$, $\{f(y_n)\}\in f^*(y)$ and $\{f(z_n)\}\in f^*(z)$. Thus by (\ref{h-1}), we have
$$T_n=\frac{\sigma_{w',\varepsilon'}(f(x_n),f(y_n))}{\sigma_{w',\varepsilon'}(f(x_n),f(z_n))}\asymp_{C(\delta)} e^{\varepsilon'[(f(x_n)|f(z_n))_{w'}-(f(x_n)|f(y_n))_{w'}]}$$
and similarly,
$$t_n=\frac{\sigma_{w,\varepsilon}(x_n,y_n)}{\sigma_{w,\varepsilon}(x_n,z_n)}\asymp_{C(\delta)} e^{\varepsilon[(x_n|z_n)_w-(x_n|y_n)_w]}.$$
It follows from Lemma \ref{z3} that there are constants $M_1\geq 1,\mu_1>0$ depending only on $M,\mu$ and $\delta$ and a function $\theta_1(t)=\max\{M_1t, t/M_1\}+\mu_1$ for $t\in \mathbb{R}$,  such that
\beq\nonumber
(f(x_n)|f(z_n))_{w'}-(f(x_n)|f(y_n))_{w'} &=& \langle f(x_n), f(z_n), f(y_n), f(w)\rangle
\\ \nonumber &\leq& \theta_1(\langle x_n,z_n,y_n,w\rangle)+\mu_1
\\ \nonumber &=& \theta_1((x_n|z_n)_w-(x_n|y_n)_w)+\mu_1.
\eeq
Therefore, combining these inequalities, we obtain
$$T_n \leq C(\delta,M,\mu,\varepsilon,\varepsilon') e^{\varepsilon'\theta_1(\varepsilon^{-1}\log t_n)}.$$
Note that $T_n\to T$ and $t_n\to t$ as $n\to \infty$. By letting $n\to \infty$, this guarantees that
$T\leq\theta(t)$ with
$$\theta(t)=C(\delta,M,\mu,\varepsilon,\varepsilon')
(t^{\frac{\varepsilon'}{\varepsilon}M_1}\vee
t^{\frac{\varepsilon'}{\varepsilon M_1}}),$$
as desired.

It remains to show that $f_*$ is relatively quasisymmetric. This can be shown by using a similar argument as above. For any sequence $\{u_n\}\in \xi$, we see from \cite[Proposition 6.3]{BS} that $\{f(u_n)\}\in \xi'=f^*(\xi)$. By (\ref{h-2}), we find that
$$S_n=\frac{\sigma_{b',\varepsilon'}(f(x_n),f(y_n))}{\sigma_{b',\varepsilon'}(f(x_n),f(z_n))}\asymp_{C(\delta)} e^{\varepsilon'[(f(x_n)|f(z_n))_{b'}-(f(x_n)|f(y_n))_{b'}]}$$
and similarly,
$$s_n=\frac{\sigma_{b,\varepsilon}(x_n,y_n)}{\sigma_{b,\varepsilon}(x_n,z_n)}\asymp_{C(\delta)} e^{\varepsilon[(x_n|z_n)_b-(x_n|y_n)_b]}.$$
Moreover, we compute

\vspace{8pt}

\noindent
$(f(x_n)|f(z_n))_{b'}-(f(x_n)|f(y_n))_{b'}$
\begin{align}
&\doteq_{C(\delta)} (f(x_n)|f(z_n))_{w'}-(f(z_n)|\xi')_{w'}-(f(x_n)|f(y_n))_{w'}+(\xi'|f(y_n))_{w'}\;\; \;(\mbox{by}\;(\ref{z0}))  \nonumber
\\ \nonumber
&\doteq_{C(\delta)} \{(f(x_n)|f(z_n))_{w'}-(f(z_n)|f(u_m))_{w'}-(f(x_n)|f(y_n))_{w'}  +(f(u_m)|f(y_n))_{w'}\}_m \\ \nonumber
&  \hspace{5cm}
 \;(\mbox{by}\;\mbox{ Lemma \ref{b-1}}) \\ \nonumber
&\doteq_{C(\delta)} \{\langle f(x_n),f(z_n),f(y_n),f(u_m) \rangle\}_m \\
\nonumber & \doteq_{C(\delta,M,\mu)} \{\theta_1(\langle x_n,z_n,y_n,u_m \rangle)\}_m\;(\mbox{by}\;   \mbox{ Lemma \ref{z3}}) \\ \nonumber
&\doteq_{C(\delta,M,\mu)} \theta_1((x_n|z_n)_b-(x_n|y_n)_b)\; \ \ (\mbox{by Lemma \ref{b-1} and (\ref{z0})})
\end{align}

Furthermore, the above arguments ensure that
$$S_n \leq C(\delta,M,\mu,\varepsilon,\varepsilon') e^{\varepsilon'\theta_1({\varepsilon}^{-1}\log s_n)}.$$
By taking $n\to \infty$, we obtain
$$\frac{\sigma_{b',\varepsilon'}(f_*(x),f_*(y))}{\sigma_{b',\varepsilon'}(f_*(x),f_*(z))}\leq \theta\Big(\frac{\sigma_{b,\varepsilon}(x,y)}{\sigma_{b,\varepsilon}(x,z)}\Big),$$
where $\theta(t)=C(\delta,M,\mu,\varepsilon,\varepsilon')
(t^{\frac{\varepsilon'}{\varepsilon}M_1}\vee
t^{\frac{\varepsilon'}{\varepsilon M_1}}).$
\epf

\begin{remark} The existence of an extension of $f^*$ on $X^*$ in Lemma \ref{zz} follows from Lemma \ref{z3}$(3)$ and \cite[Corollary 4.3.3]{BuSc} or \cite[Proposition 6.3]{BS}.
\end{remark}
In \cite{BuSc}, Buyalo and Schroeder introduced a special kind of hyperbolic fillings, which was known as {\it hyperbolic approximation}, for complete metric spaces. Many  applications are presented with the aid of this useful construction. For instance, they proved bilipschitz and quasisymmetric extension results on Gromov hyperbolic spaces from the boundaries at infinity to their interiors (see \cite[Chapter 7]{BuSc}), Bonk-Schramm Embedding Theorem (see \cite{BS}) and  Assouad Embedding Theorem (see \cite[Chapter 8]{BuSc}).

Next, using the hyperbolic approximation structures, we extend Lemma \ref{zz} to a more general context as follows.

\begin{lem}\label{w-0}
Suppose that $X$ is BuSc-visual and $Y$ is rough geodesic, and that $f:X\to Y$ is a rough quasi-isometry which induces a bijective map $f^*: \partial^* X\to \partial^* Y$.
Let $w'=f(w)$ and $\xi'=f^*(\xi)$. Then we have the following:
\begin{enumerate}
\item $f^*:(X^*,\sigma_{w,\varepsilon})\to (Y^*,\sigma_{w',\varepsilon'})$ is relatively quasisymmetric;
\item $f_*:(X^*\setminus\{\xi\},\sigma_{b,\varepsilon})\to (Y^*\setminus\{\xi'\},\sigma_{b',\varepsilon'})$ is relatively quasisymmetric.
\end{enumerate}
\end{lem}
\bpf To prove the first assertion, we consider the hyperbolic approximations of bounded complete metric spaces $(\partial^* X,\sigma_{w,\varepsilon})$ and $(\partial^* Y,\sigma_{w',\varepsilon'})$, and let
$$Hyp_r(\partial^* X,\sigma_{w,\varepsilon})=:\widetilde{X_r}\;\;\mbox{ and }\;\;\; Hyp_{r'}(\partial^* Y,\sigma_{w',\varepsilon'})=:\widetilde{Y_{r'}}$$
for some parameters $0<r,r'<1$, respectively. For our discussion, we recall certain useful properties of $\widetilde{X_r}$ and $\widetilde{Y_{r'}}$. According to \cite[Theorem $3$]{JJ}, it follows that both of them are BuSc-visual geodesic Gromov hyperbolic spaces. There is a canonical bilipschitz identification between $(\partial^* \widetilde{X_r},\sigma_{o,-\ln r})$ and $(\partial^* X,\sigma_{w,\varepsilon})$. Similarly, the canonical identification between $(\partial^* \widetilde{Y_{r'}},\sigma_{o',-\ln r'})$ and $(\partial^* Y,\sigma_{w',\varepsilon'})$ is bilipschitz.

Applying the assumption and \cite[Theorem $7.1.2$]{BuSc} or \cite[Theorem $1$]{JJ}, the above facts imply that there are roughly isometric maps
$\varphi: X\to \widetilde{X_r}$ and $\psi: \widetilde{Y_{r'}}\to Y$
such that $\varphi^*|_{\partial^*X}=\id_{\partial^*X}$ and $\psi^*|_{\partial^*Y}=\id_{\partial^*Y}$.

\begin{figure}[htbp]
\begin{center}
\includegraphics{Lemma-4-8.1}
\end{center}
\caption{} \label{2}
\end{figure}

Furthermore, since $f:X\to Y$ is a rough quasi-isometry, there is a rough inverse $g:Y\to X$ of $f$, which is also a rough quasi-isometry, such that $f\circ g\doteq \id_Y$ and $g\circ f\doteq \id_X$ (cf. \cite{Vai-0}). Also we see from the assumption that $g^*|_{\partial^*Y}=(f^*|_{\partial^*X})^{-1}$ is a bijection. Hence, we obtain a roughly quasi-isometric map
$$\widetilde{g}=\varphi\circ g \circ \psi: \widetilde{Y_{r'}}\to \widetilde{X_r}$$
such that $\widetilde{g}^*|_{\partial^*\widetilde{Y_{r'}}}=(\varphi\circ g \circ \psi)^*|_{\partial^*\widetilde{Y_{r'}}}=(\varphi^*\circ g^* \circ \psi^*)|_{\partial^*\widetilde{Y_{r'}}}$ is a bijection. It follows from \cite[Lemma $7.3.12$]{BuSc}  that $\widetilde{g}$ is a rough quasi-isometry between these two hyperbolic approximation spaces. Thus, we find that there is a rough inverse
$\widetilde{f}:\widetilde{X_r}\to \widetilde{Y_{r'}}$
of $\widetilde{g}$, which is a rough quasi-isometry, such that $\widetilde{f}\circ \widetilde{g} \doteq \id_{\widetilde{Y_{r'}}}$ and $\widetilde{g}\circ \widetilde{f}\doteq \id_{\widetilde{X_r}}$. We claim that
\be\label{b-2}\psi\circ \widetilde{f}\circ \varphi \doteq f.\ee
Indeed, we note that
$$(\varphi\circ g \circ \psi)\circ \widetilde{f}\circ \varphi=\widetilde{g}\circ \widetilde{f}\circ \varphi \doteq \id_X \circ \varphi=\varphi.$$
Since $\varphi$ is a rough isometry, we have $g \circ \psi \circ \widetilde{f}\circ \varphi\doteq  \id_X$ and therefore,
\beq\nonumber
f=f\circ \id_X &\doteq& f\circ (g \circ \psi \circ \widetilde{f}\circ \varphi)
\\ \nonumber&=& (f\circ g)\circ \psi \circ \widetilde{f}\circ \varphi
\\ \nonumber&\doteq& \id_Y\circ \psi\circ \widetilde{f}\circ \varphi=\psi\circ \widetilde{f}\circ \varphi.
\eeq
Hence, (\ref{b-2}) is true.

Now we are ready to check the first assertion. By (\ref{b-2}), to prove the relative quasisymmetry of $f^*$,
we only need to show that
$$(\psi\circ \widetilde{f}\circ \varphi )^*=\psi^*\circ \widetilde{f}^*\circ \varphi ^*:(X^*,\sigma_{w,\varepsilon})\to (Y^*,\sigma_{\widetilde{w'},\varepsilon'})$$
is relatively quasisymmetric, where $\widetilde{w'}= \psi \circ \widetilde{f}\circ \varphi(w)$. It thus follows from Lemma \ref{z5} and the fact that $\psi$, $\varphi$ are roughly isometric maps that $\psi^*$ and $\varphi ^*$ are relatively snowflake. Combining Lemma \ref{zz} and the fact that $\widetilde{f}$ is a rough quasi-isometry between geodesic Gromov hyperbolic spaces, we obtain that $\widetilde{f}^*$ is relatively quasisymmetric. Therefore, the composition $\psi^*\circ \widetilde{f}^*\circ \varphi ^*$ is also relatively quasisymmetric, as required.

It remains to show the second assertion. Lemma \ref{z1} guarantees that the identity map
$\id_{X^*}:(X^*\setminus \{\xi\},\sigma_{w,\varepsilon})\to (X^*\setminus\{\xi\},\sigma_{b,\varepsilon})$
is quasim\"{o}bius, and $(X^*\setminus \{\xi\},\sigma_{b,\varepsilon})$ is unbounded with $\sigma_{b,\varepsilon}(\xi,\zeta)=\infty$ for all $\xi\neq \zeta\in\partial^*X$. Similarly, we observe that the identity map
$\id_{Y^*}:(Y^*\setminus \{\xi'\},\sigma_{w',\varepsilon'})\to (Y^*\setminus\{\xi'\},\sigma_{b',\varepsilon'})$
is quasim\"{o}bius with $\sigma_{b',\varepsilon'}(\xi',f^*(\zeta))=\infty$.

Therefore, it follows that
$$f_*=\id_{Y^*}\circ f^* \circ (\id_{X^*})^{-1}:(X^*\setminus\{\xi\},\sigma_{b,\varepsilon})\to (Y^*\setminus\{\xi'\},\sigma_{b',\varepsilon'})$$
is relatively quasim\"{o}bius, because the compositions of quasim\"obius maps and relatively quasisymmetric maps are relatively quasim\"{o}bius. Hence, Lemma \ref{z-1} ensures that $f_*$ is in fact a relatively quasisymmetric map. We can complete the proof.
\epf

As an immediate consequence of Lemma \ref{w-0}, we investigate the invariance of visual property of Gromov hyperbolic spaces under roughly quasi-isometric maps. For later use, we need the following result which may be of independent interest.

\begin{thm}\label{zz0} Let $X$ be a BuSc-visual $\delta$-hyperbolic space and $Y$ be a roughly geodesic $\delta$-hyperbolic space. If $f:X\to Y$ is a rough quasi-isometry which induces a bijective map $f^*: \partial^* X\to \partial^* Y$, then $Y$ is BoSc-visual.
\end{thm}
\bpf Assume that $X$ is $K$-BuSc-visual with respect to $w\in X$ for some constant $K\geq 0$. To show that $Y$ is BoSc-visual, by Lemma \ref{g-5}, we only need to prove that $Y$ is BuSc-visual with respect to $w'=f(w)$. For each $y\in Y$, there is an $x\in X$ and a number $C\geq 0$ such that
$d'(y,f(x))\leq C,$
because $f$ is a rough quasi-isometry. Since  $X$ is BuSc-visual, there is a $\xi\in\partial^*X$ such that
$d(w,x)\leq (x|\xi)_w+K.$
For any sequence $\{a_n\}\in \xi$, we have $\{f(a_n)\}\in f^*(\xi)=\xi'$. Without loss of generality, we may assume for all $n$ sufficiently large that, $a_n\not\in\{x,w\}$ and $f(a_n)\not\in\{f(x),f(w)\}$. By Lemma \ref{b-1}, it follows that there is a constant $K_0=K_0(K,\delta)$ such that
$(w|a_n)_x=d(x,w)-(x|a_n)_w\leq K_0.$

On the other hand, Lemma \ref{w-0} implies that there is a homeomorphism $\theta:[0,\infty)\to [0,\infty)$ such that the induced map $f^*:(X^*,\sigma_{x,\varepsilon})\to (Y^*,\sigma_{f(x),\varepsilon'})$ is relatively $\theta$-quasisymmetric, where $\sigma_{x,\varepsilon}$ and $\sigma_{f(x),\varepsilon'}$ are hyperbolic metrics of $X$ and $Y$, respectively. Thus, by (\ref{h-1}), we obtain
\beq\nonumber
\frac{\sigma_{f(x),\varepsilon'}(f(x),f(a_n))}{\sigma_{f(x),\varepsilon'}(f(w),f(a_n))} &\asymp_C& e^{\varepsilon'[(f(w)|f(a_n))_{f(x)}- (f(x)|f(a_n))_{f(x)}]}=e^{\varepsilon'(f(w)|f(a_n))_{f(x)}}
\\ \nonumber&\leq& \theta\Big(\frac{\sigma_{x,\varepsilon}(x,a_n)}{\sigma_{x,\varepsilon}(w,a_n)}\Big) \asymp_C \theta(e^{\varepsilon (w|a_n)_x})
\\ \nonumber&\leq& \theta(e^{\varepsilon K_0}),
\eeq
which implies that
$(f(w)|f(a_n))_{f(x)}\leq K_1$
for some constant $K_1=K_1(\varepsilon,\varepsilon',\delta,\theta,K_0)\geq 0$. Therefore, by \cite[Lemma $2.8(5)$]{Vai-0}, we have
\beq\nonumber
d'(w',y)&\leq& d'(w',f(x))+d'(f(x),y)
\\ \nonumber&=& (f(w)|f(a_n))_{f(x)}+(f(x)|f(a_n))_{f(w)}+d'(f(x),y)
\\ \nonumber&\leq& (f(w)|f(a_n))_{f(x)}+(y|f(a_n))_{f(w)}+2d'(f(x),y)
\\ \nonumber&\leq& (y|f(a_n))_{f(w)}+K_1+2C.
\eeq
Consequently, we obtain from Lemma \ref{b-1} that there is a constant $K'=K'(K_1,C,\delta)$ such that
$d'(y,f(w))\leq (f^*(\xi)|y)_{f(w)}+K'.$
This completes the proof.
 \epf

We conclude this section by considering quasisymmetric extension properties of Gromov hyperbolic spaces. Lemma \ref{z6} below is related to \cite[Theorem 1]{JJ}.

\begin{lem}\label{z6}
Let $X$ be a BuSc-visual Gromov hyperbolic space and $Y$ be a roughly geodesic Gromov hyperbolic space. If there is a power quasisymmetric homeomorphism $\phi:(\partial^*X,\sigma_{w,\varepsilon})\to (\partial^* Y,\sigma_{w',\varepsilon'}),$
then $\phi$ induces a roughly quasi-isometric map $f:X\to Y$ such that $f^*:(X^*,\sigma_{w,\varepsilon})\to (Y^*,\sigma_{f(w),\varepsilon'})$ is relatively quasisymmetric with $f^*|_{\partial^*X}=\phi$.
\end{lem}
\bpf Note that our arguments are similar to the proof of Lemma \ref{w-0}. For completeness we give the details.

Firstly, we consider the hyperbolic approximations of bounded complete metric spaces $(\partial^* X,\sigma_{w,\varepsilon})$ and $(\partial^* Y,\sigma_{w',\varepsilon'})$, denoted by
$$Hyp_r(\partial^* X,\sigma_{w,\varepsilon})=:\widetilde{X_r}\;\;\mbox{ and }\;\;\; Hyp_{r'}(\partial^* Y,\sigma_{w',\varepsilon'})=:\widetilde{Y_{r'}}$$
for some parameters $0<r,r'<1$, respectively. We also record some useful properties of $\widetilde{X_r}$ and $\widetilde{Y_{r'}}$. By \cite[Theorem $3$]{JJ}, both of them are BuSc-visual geodesic Gromov hyperbolic spaces. Moreover, there is a canonical bilipschitz identification between $(\partial^* \widetilde{X_r},\sigma_{o,-\ln r})$ and $(\partial^* X,\sigma_{w,\varepsilon})$. Similarly, there is a canonical bilipschitz identification between $(\partial^* \widetilde{Y_{r'}},\sigma_{o',-\ln r'})$ and $(\partial^* Y,\sigma_{w',\varepsilon'})$.

\begin{figure}[htbp]
\begin{center}
\includegraphics{Lemma-4-10.1}
\end{center}
\caption{} \label{3}
\end{figure}

On the one hand, we see from the above facts that $\phi$ induces a power quasisymmetric homeomorphism
$$\widetilde{\phi}=\id_{\partial^*Y}\circ \phi\circ (\id_{\partial^*X})^{-1}: (\partial^* \widetilde{X_r},\sigma_{o,-\ln r})\to (\partial^* \widetilde{Y_{r'}},\sigma_{o',-\ln r'}),$$
because the compositions of power quasisymmetric maps and bilipschitz maps are power quasisymmetric. According to \cite[Theorem $1$]{JJ}, $\widetilde{\phi}$ induces a roughly quasi-isometric map
$\widetilde{f}:\widetilde{X_r}\to \widetilde{Y_{r'}} $
with $\widetilde{f}^*|_{\partial^* \widetilde{X_r}}=\widetilde{\phi}$. Thus by Lemma \ref{zz}, it follows that $\widetilde{f}$ induces a relatively quasisymmetric map
$\widetilde{f}^*:(\widetilde{X_r}^*,\sigma_{\varphi(w),\widetilde{\varepsilon}})\to (\widetilde{Y_{r'}}^*,\sigma_{\widetilde{f}\circ \varphi(w),\widetilde{\varepsilon'}}).$

On the other hand, since $X$ is BuSc-visual $\delta$-hyperbolic and $Y$ is rough geodesic $\delta$-hyperbolic, using the properties of hyperbolic approximations  and \cite[Theorem $1$]{JJ}, we obtain that there are roughly isometric maps
$\varphi: X\to \widetilde{X_r}$ and $\psi: \widetilde{Y_{r'}}\to Y,$
with $\varphi^*|_{\partial^*X}=\id_{\partial^*X}$ and $\psi^*|_{\partial^*Y}=\id_{\partial^*Y}$, respectively. Now Lemma \ref{z5} ensures that
$$\varphi^*: (X^*,\sigma_{w,\varepsilon})\to (\widetilde{X_r}^*,\sigma_{\varphi(w),\widetilde{\varepsilon}})\;\;\;\;\;\; \mbox{and}\;\;\;\;\;\; \psi ^*:(\widetilde{Y_{r'}}^*,\sigma_{\widetilde{f}\circ \varphi(w),\widetilde{\varepsilon'}})\to (Y^*,\sigma_{\widetilde{w'},\varepsilon'})$$
are relatively snowflake, where $\widetilde{w'}=\psi\circ \widetilde{f}\circ \varphi(w)=f(w)$.

Let
$f=\psi\circ \widetilde{f} \circ \varphi:X\to Y.$
One easily sees that $f$ is roughly quasi-isometric. Moreover, we find that the induced mapping
$f=\psi^*\circ\widetilde{f}^*\circ \varphi^*:(X^*,\sigma_{w,\varepsilon})\to (Y^*,\sigma_{f(w),\varepsilon'}).$
is relatively quasisymmetric because the compositions of relatively snowflake maps and relatively quasisymmetric maps are relatively quasisymmetric.

It thus left to show that $f^*|_{\partial^*X}=\phi$. Indeed, we have
$$f^*|_{\partial^*X}=\psi^*\circ\widetilde{f}^*\circ \varphi^*|_{\partial^*X}=\psi^*\circ \id_{\partial^*Y}\circ \phi\circ (\id_{\partial^*X})^{-1}\circ \varphi^*|_{\partial^*X}=\phi,$$
as desired. Hence the lemma follows.
\epf

\section{Hyperbolic type metrics on uniform spaces}
This section focuses on two hyperbolic type metrics (in the sense of Gromov): the hyperbolization metric $h$ and the quasihyperbolic metric $k$. Throughout this section, we do not assume that the spaces are locally compact unless otherwise specified.

\subsection{Hyperbolization of incomplete spaces}

In \cite{Z1}, Ibragimov introduced a metric that hyperbolizes $($in the sense of Gromov$)$ all incomplete and locally complete metric spaces, which is roughly quasi-isometric to the quasihyperbolic metric on uniform metric spaces. Here, a metric space $(X, d)$ is called {\it locally complete} if each point has a complete neighborhood. Note that locally compact spaces are locally complete.

Let $(X,d)$ be an incomplete and locally complete metric space. We define the {\it hyperbolization  metric} $h$ of $X$ as
\be\label{h} h(x,y)=2\log\frac{d(x,y)+ d(x)\vee d(y)}{\sqrt{d(x)d(y)}},\ee
for all $x,y\in X.$ Next, we recall certain useful results about this metric.

\begin{thm}\label{zz2}  Let $(X,d)$ be an incomplete and locally complete metric space.
\begin{enumerate}
\item $($\cite[Theorem $2.1$]{Z1}$)$ The space $(X,h)$ is $\delta$-hyperbolic with $\delta=\log 4$.
\item $($\cite[Theorem $3.1$]{Z1}$)$ For all $x,y\in X$, $j(x,y)\leq h(x,y)\leq 2j(x,y)+2\log 2$.
\item $($\cite[Theorem $3.1$]{Z2}$)$ The boundary at infinity $\partial_h^* X$of $(X,h)$ is naturally equivalent to its metric boundary $\partial X~ (\partial X \cup \{\infty\}$ if $X$ is unbounded$)$.
\end{enumerate}
\end{thm}

In the following, we show that the hyperbolization space $(X,h)$ is BuSc-visual.  We first consider the bounded case.

\begin{lem}\label{g-3}  Let $(X,d)$ be a bounded incomplete and locally complete metric space. Then $(X,h)$ is a BuSc-visual hyperbolic space.
\end{lem}
\bpf By Theorem \ref{zz2}, $(X,h)$ is a $\delta$-hyperbolic space with $\delta=\log4$. We only need to prove that there is a positive constant $K$ and a fixed point $w$ in $X$ such that $(X,h)$ is $K$-BuSc-visual with respect to $w$.

To this end, without loss of generality we assume that there is a point $w\in X$ such that
$d(w)=\sup\{d(x):x\in X\}.$
Fix $y\in X$. Choose a point $y'\in \partial X$ satisfying $d(y,y')\leq 2d(y)$.
Then take a sequence of points $\{y_n\}$ in $X$ such that $d(y_n,y')<d(y)/2$ and $d(y_n, y')\to 0$, as $n\rightarrow\infty$. By Theorem \ref{zz2}, we may regard the boundary at infinity $\partial_h^* X$ of $(X,h)$ and the metric boundary $\partial X$ of $(X,d)$ as a set for the simplicity of notations.

Next, we shall find an upper bound for the quantity
$h(w,y)-(y|y_n)_w=(w|y_n)_y.$
Since $y'\in \partial X$ and $d(y_n, y')\to 0$ as $n\rightarrow\infty$, we have $d(y_n)\to 0$ as $n\rightarrow\infty$. Thus we may assume that for all $n$, $d(y_n)\leq d(y)\leq d(w)$. By (\ref{h}), an elementary computation guarantees that
\begin{eqnarray*} 2(w|y_n)_y &=& h(w,y)+h(y_n,y)-h(w,y_n)
\\ \nonumber&=& 2\log\frac{[d(w,y)+d(x)\vee d(y)][d(y_n,y)+d(y)]}{d(y)[d(y_n,w)+d(w)]}
.\end{eqnarray*}
By the choices of $w$, $y'$ and $y_n$, we have
$$d(y_n,y)+d(y)\leq  d(y_n,y')+d(y',y)+d(y)\leq d(y)/2+2d(y)+d(y)= 7d(y)/2,$$
and
$$d(w,y)+d(x)\vee d(y) \leq  d(w,y_n)+d(y_n,y)+d(w)\leq 4(d(y_n,w)+d(w)).$$
It follows from the above two estimates that
$(w|y_n)_y \leq \log 14.$
This, together with Lemma \ref{b-1}, ensures that
$$h(w,y)-(y|y')_w\leq {(w|y_n)_y}_n+2\delta\leq \log 14+4\log 2
$$
The proof of Lemma \ref{g-3} is complete.
 \epf

For the unbounded case, our result is as follows.

\begin{lem}\label{z-2}  Let $(X,d)$ be an unbounded incomplete and locally complete metric space. Then $(X,h)$ is BuSc-visual hyperbolic with respect to any point of $X$.
\end{lem}
\bpf For a given point $w\in X$, we show that $(X,h)$ is BuSc-visual with respect to $w$. For each $y\in X$ and any sequence $\{y_n\}\subset X$, a direct computation gives that
\beq\label{b-4} 2(w|y_n)_y&=& h(w,y)+h(y_n,y)-h(w,y_n)
\\ \nonumber&=& 2\log\frac{[d(w,y)+d(w)\vee d(y)][d(y_n,y)+d(y_n)\vee d(y)]}{d(y)[d(y_n,w)+d(y_n)\vee d(y)]}
.\eeq

By Theorem \ref{zz2}, it follows that there is a natural identification
$\varphi:\partial_h^* X\to\partial X\cup\{\infty\}$
between the boundary at infinity $\partial_h^* X$ of $(X,h)$ and the one-point extension $\partial X\cup\{\infty\}$ of $\partial X$ with $\varphi(\xi)=\infty$ for some $\xi\in \partial_h^* X$. We consider two cases in the following.

Suppose first that $3d(y)\leq (d(w)\vee d(w,y)).$ For each $\epsilon\in(0,1)$, choose a point $y'\in \partial X$ satisfying $d(y,y')<(1+\epsilon)d(y)$. Take a sequence of points $\{y_n\}$ in $X$ such that $d(y_n,y')<\epsilon d(y)$ and $d(y_n, y')\to 0$ as $n\rightarrow\infty$. Next, we shall find an upper bound for $h(w,y)-(y|y_n)_w=(w|y_n)_y$.

By the choices of $y'$ and $y_n$, we have
$$d(y_n,y)\leq d(y,y')+d(y_n,y')<(1+2)\epsilon d(y),
$$
$$d(y_n)\vee d(y)\leq  d(y_n,y')+d(y',y)+d(y)\leq 2(1+\epsilon)d(y),$$
and
\beqq
d(w,y_n)+d(y_n)\vee d(w)&\geq&  d(w,y)-d(y',y)-d(y_n,y')+d(w)
\\ \nonumber &\geq& d(w,y)+d(w)-(1+2\epsilon)d(y)
\\ \nonumber &\geq& \frac{2}{3}(1-\epsilon)(d(w,y)+d(w))
\\ \nonumber &\geq&  \frac{1}{3}(1-\epsilon)(d(w,y)+ d(w)\vee d(y)  ).
\eeqq
Using these inequalities and (\ref{b-4}), we obtain
$$(w|y_n)_y \leq \log\frac{12(1+\epsilon)}{1-\epsilon}.$$
By letting $\epsilon\to 0$, Lemma \ref{b-1} ensures that
$$h(w,y)-(y|y')_w\leq \{(w|y_n)_y\}_n+2\delta\leq \log3+6\log 2=K ~{\rm (say)}.
$$
Hence we are done.

Next, we consider the case that $3d(y)>  d(w)\vee d(w,y).$ Thus we have
\be\label{b-1.1} d(w,y)+  d(w)\vee d(y)  <4d(y).\ee
Because $X$ is unbounded, there is a sequence of points $\{z_n\}\subset X$ with $d(z_n,y)\to\infty$ as $n\to\infty$. Thus, we know that $\{z_n\}\in\xi$.
By (\ref{b-4}) and (\ref{b-1.1}), we compute from the triangle inequality that for all $n\in \mathbb{N}$,
$$(w|z_n)_y\leq \log \frac{4[d(z_n,y)+d(y)+d(z_n)]}{d(z_n,y)-d(y,w)+d(z_n)}.$$
Moreover, by taking $n\to \infty$, we see from Lemma \ref{b-1} that
$$h(w,y)-(y|\xi)_w\leq \{(w|z_n)_y\}_n+2\delta\leq 6\log 2\leq K.$$
The proof is complete.
\epf

Next, for an incomplete and locally complete metric space $(X,d)$, we prove that there is a naturally relatively quasim\"{o}bius identification between the Gromov closure of the hyperbolization space $(X,h)$ and $\overline{X}'$. If, in addition, $X$ is uniform, the natural map is relatively quasisymmetric, quantitatively. These results are useful to investigate the relationships between the Gromov boundary associated with the quasihyperbolic metric and the metric boundary in the sequel.

\begin{lem}\label{g-1} Let $(X,d)$ be an incomplete and locally complete metric space with $w\in X$, and let $0<\varepsilon\leq 1/2$ and  $0<\varepsilon'\leq 1/44$.
\begin{enumerate}
\item\label{sun-4} The natural bijective map $\varphi: (X_h^*,h_{w,\varepsilon})\rightarrow (\overline{X}',d)$ exists and is relatively $\theta$-quasim\"obius with $\theta(t)=16t^{\varepsilon}$. Note that $\overline{X}'$ is defined as in Remark \ref{r-1}.
\item\label{sun-5} If $X$ is bounded and if $w\in X$ with $d(w)=\sup\{d(x):x\in X\}$, then $\varphi$ is relatively $\eta$-quasisymmetric with $\eta(t)=\frac{12\diam(X)}{d(w)}t^{\varepsilon}$. In particular, if, in addition, $X$ is $A$-uniform, then the control function takes the form $\eta(t)=24At^{\varepsilon}$.
\item\label{sun-6} If $X$ is unbounded, then the natural map $\varphi: (X_h^*\setminus\{\xi\},h_{b,\varepsilon'})\rightarrow (\overline{X},d)$ is relatively $\eta'$-quasisymmetric with $\eta'(t)=64t^{\varepsilon'}$, where $b=b_{\xi,w}\in \mathcal{B}_h(\xi)$ is a Busemann function based at $\xi$ with $\varphi(\xi)=\infty$.
    \end{enumerate}
\end{lem}

\bpf
Firstly, Theorem \ref{zz2} implies that $(X,h)$ is $\delta$-hyperbolic with $\delta=\log4$ and there exists a natural bijective map $\varphi: (X_h^*,h_{w,\varepsilon})\rightarrow (\overline{X}',d)$. For the purpose of simplicity in notation, we consider the sets $\partial_h^*X$ and $\partial X$ to be the same. Therefore, we know that a sequence of points $\{x_n\}$ in $X$ $d$-converges to a point $x\in \partial X'$ if and only if $\{x_n\}$ is a Gromov sequence with $\{x_n\}\in x\in \partial_h^*X$.

To verify \eqref{sun-4}, we only need to show that $\varphi$ is relatively $\theta$-quasim\"obius. For all $x,$ $y\in X$, a direct computation gives that
\beq\label{g-0.4}\;\;\;\;\;\;\; (x|y)_w &=& \frac{1}{2}\big(h(x,w)+h(y,w)-h(x,y)\big)
\\ \nonumber  &=& \log\frac{[d(w,x)+ d(w)\vee d(x) ][d(w,y)+d(w)\vee d(y)]}
{d(w)[d(x,y)+d(x)\vee d(y)]}.\eeq

Let $x,y,z,v$ be four distinct  points in $\overline{X}$ with $x,v\in \partial_h^*X$, $\{x_n\}\subset X$ and $\{y_n\}\subset X$ (if $y\in X$, let $y_n=y$ for all $n\geq 1$) such that $d(x_n,x)\rightarrow 0$ and $d(y_n,y)\rightarrow 0$ as $n\rightarrow\infty$, respectively. By (\ref{g-0.4}), we have
$$(x_n|y_n)_w \rightarrow \log \frac{[d(w,x)+d(w)][d(w,y)+d(w)\vee d(y)]}{d(w)(d(x,y)+d(y))} ~\mbox{ as $n\rightarrow\infty$}.$$
This leads to
\be\label{z-0}\rho_{w,\varepsilon}(x,y)=e^{-\varepsilon(x|y)_w}
=\bigg(\frac{d(w)(d(x,y)+d(y))}{(d(w,x)+
d(w))(d(w,y)+d(w)\vee d(y))}\bigg)^{\varepsilon},\ee
where $d(y)=0$ if $y\in \partial X$. Note that a similar argument as above gives the corresponding formula for the quantities $\rho_{w,\varepsilon}(x,z)$, $\rho_{w,\varepsilon}(y,v)$ and  $\rho_{w,\varepsilon}(z,v)$. Using these formula and (\ref{h-1}), we see that the cross ratio of $E=\{x,y,z,v\}$ with respect to the hyperbolic metric $h_{w,\varepsilon}$ satisfies
\begin{eqnarray*}\tau(E,h_{w,\varepsilon}) &=& \frac{h_{w,\varepsilon}(x,z)h_{w,\varepsilon}(y,v)}{h_{w,\varepsilon}(x,y)h_{w,\varepsilon}(z,v)}
\leq \frac{4\rho_{w,\varepsilon}(x,z)\rho_{w,\varepsilon}(y,v)}{\rho_{w,\varepsilon}(x,y)\rho_{w,\varepsilon}(z,v)}
\\ \nonumber&=& 4\bigg(\frac{(d(x,z)+d(z))(d(y,v)+d(y))}{(d(x,y)+d(y))(d(z,v)+d(z))}\bigg)^{\varepsilon}
\leq  16\tau(E,d)^{\varepsilon},\end{eqnarray*}
where $\tau(E,d)$ denotes the cross ratio of $E=\{x,y,z,v\}$ with respect to the metric $d$. Therefore, it follows from \cite[Remark 6.21]{Vai-4} that $\varphi$ is relatively $\theta$-quasim\"obius with $\theta(t)=16t^{\varepsilon}$,  which shows that \eqref{sun-4} is true.

Next, we check \eqref{sun-5}. Suppose that $X$ is bounded. To prove that the natural map $\varphi$ is relatively $\eta$-quasisymmetric, let $x,y,z$ be three distinct arbitrary points in $X_h^*$ with $x\in \partial_h^*X$. By the choice of $w$, it follows from a similar argument as in the case of the proof of (\ref{z-0}) that
\begin{eqnarray*} \frac{h_{w,\varepsilon}(x,y)}{h_{w,\varepsilon}(x,z)} &\leq& \frac{2\rho_{w,\varepsilon}(x,y)}{\rho_{w,\varepsilon}(x,z)}
\\ \nonumber&=& 2\bigg(\frac{d(x,y)+d(y)}{d(x,z)+d(z)}\bigg)^{\varepsilon}
\bigg(\frac{d(w,z)+ d(w)\vee d(z)}{d(w,y)+d(w)\vee d(y)}\bigg)^{\varepsilon}
\\ \nonumber&\leq& \frac{12\diam(X)}{d(w)}\bigg(\frac{d(x,y)}{d(x,z)}\bigg)^{\varepsilon}.\end{eqnarray*}
This ensures that $\varphi$ is relatively $\eta$-quasisymmetric with $\eta(t)=\frac{12\diam(X)}{d(w)}t^{\varepsilon}$ by invoking \cite[Lemma 6.17]{Vai-4}.

Assume further that $X$ is $A$-uniform. For all $p,q\in X$, we know that there is an $A$-uniform arc $\alpha$ connecting $p$ and $q$. Let $x_0\in \alpha$ be a point such that
$\ell(\alpha[p,x_0])=\ell(\alpha[q,x_0])$. By the uniformity of $\alpha$ and by the choice of $w$, we get
$$d(p,q)\leq \ell(\alpha)\leq 2Ad(x_0)\leq 2Ad(w),
$$
which implies that
$\diam(X)\leq 2A d(w).$
Therefore, we may take $\eta(t)=24At^{\varepsilon}$ and hence, \eqref{sun-5} holds true.

All that remains to be proved is \eqref{sun-6}. Suppose that $X$ is unbounded.
Let $x,y,z$ be three distinct points in $X_h^*\setminus\{\xi\}$ with $x\in \partial_h^*X\setminus\{\xi\}$. For any sequences $\{x_n\}\in x$ and $\{u_n\}\in\xi$,
we thus have $d(x_n,x)\to0$ and $d(u_n,y)\to \infty$ as $n\to \infty$. Without loss of generality we may assume that $y,z\in X$ and $d(x_n)\leq  d(y)\wedge d(z)$ for all $n$ sufficiently large. It follows that
\beq\nonumber e^{(x_n|z)_b-(x_n|y)_b} &\leq& e^{20\delta}e^{(x_n|z)_w-(x_n|y)_w-(\xi|z)_w-(\xi|y)_w}\;\;\;\;\;\;\;\;\;\;\;\;\;\;\;\;\;\;\;\;\;\;\;\;\;\;\;\;
(\mbox{by}\;(\ref{z0}))
\\ \nonumber  &\leq& e^{24\delta}e^{(x_n|z)_w-(x_n|y)_w-(u_n|z)_w-(u_n|y)_w}\;\;\;\;\;\;\;\;\;\;\;\;\;\;\;\;\;\;\;\;\;\;\;\;
(\mbox{by\;Lemma \ref{b-1}})
\\ \nonumber  &=& e^{24\delta} \frac{d(x_n,y)+d(y)}{d(x_n,z)+d(z)} \frac{d(u_n,z)+d(u_n)\vee d(z)}{d(u_n,y)+d(u_n)\vee d(y)}  \;\;\;\;\;\;\;\;(\mbox{by}\;(\ref{g-0.4}))
\\ \nonumber  &\leq& e^{24\delta} \frac{2d(x_n,y)+d(x_n)}{d(x_n,z)} \frac{2d(u_n,z)+d(u_n)}{d(u_n,y)+d(u_n)}.
\eeq
Note that $0<\varepsilon'\leq 1/44$ and $\delta=\log4$. We further obtain
\beq\nonumber \frac{h_{b,\varepsilon'}(x,y)}{h_{b,\varepsilon'}(x,z)} &\leq& 2e^{\varepsilon'[(x|z)_b-(x|y)_b]}\;\;\;\;\;\;\;\;\;\;\;\;\;\;\;\;\;\;\;\;\;\;\;\;\;\;\;\;\;(\mbox{by}\;(\ref{h-2}))
\\ \nonumber  &\leq& 2e^{22\delta\varepsilon'}\{e^{\varepsilon'[(x_n|z)_b-(x_n|y)_b]}\}_n\;\;\;\;\;\;\;\;\;\;\;\;\;(\mbox{by \cite[Lemma 3.2.4]{BuSc}})
\\ \nonumber  &\leq& 2e^{44\delta\varepsilon'} \bigg(\frac{4d(x,y)}{d(x,z)}\bigg)^{\varepsilon'}\;\;\;\;\;\;\;\;\;\;\;\;\;\;\;\;\;\;\;\;\;\;\;(\mbox{by the above inequality})
\\ \nonumber  &\leq& 32 \bigg(\frac{d(x,y)}{d(x,z)}\bigg)^{\varepsilon'},
\eeq
as needed. This completes the proof.
\epf

\subsection{Quasihyperbolization of uniform spaces}
In this part, using the results about the hyperbolization metric $h$ established in the preceding subsection, we provide a different method to study the Gromov hyperbolicity of uniform spaces with respect to the quasihyperbolic metric.

\begin{lem}\label{g-0}  If $(X,d)$ is an $A$-uniform metric space, then $(X,k)$ is $\delta$-hyperbolic with $\delta=\delta(A)$.
\end{lem}

\bpf Firstly, we observe from Theorem \ref{zz2} that $(X,h)$ is $\delta_0$-hyperbolic with $\delta_0=\log4$. It follows from \cite[Theorem 4.1]{BS} that $(X,h)$ can embed isometrically into a geodesic Gromov hyperbolic metric space $\widetilde{X}$. Moreover, according to  Theorems \ref{zz1} and \ref{zz2}, we find that the identity map $(X,k)\to (X,h)$ is $(M,\mu)$-roughly quasi-isometric with the constants $M,\mu$ depending only on $A$. So $(X,k)$ is roughly quasi-isometric embedding into the space $\widetilde{X}$. It  follows from \cite[Theorem 3.18]{Vai-0} that $(X,k)$ is $\delta$-hyperbolic with $\delta=\delta(\lambda,\mu,\delta_0)=\delta(A)$, because evidently $(X,k)$ is a length space. Note that in \cite[Theorem 3.18]{Vai-0}, the roughly quasi-isometric map do not need to satisfy the roughly surjective condition.
\epf

A similar argument as in the proof of \cite[Theorem 1.3]{Z-h} gives the following result.

\begin{lem}\label{zz-08}  If $(X,d)$ is an $A$-uniform metric space, then the natural map $\varphi: X_k^*\rightarrow \overline{X}'$ exists and is a bijection.
\end{lem}

Next, we show that bounded uniform spaces are visual which is a generalization of \cite[Theorem 3.6]{BHK}. Our method of proof is different.

\begin{lem}\label{h-03} If $(X,d)$ is a bounded $A$-uniform metric space, then $(X,k)$ is a BuSc-visual and BoS-visual hyperbolic space.
\end{lem}
\bpf It follows from Lemmas \ref{g-1} and \ref{zz-08} that the identity map $\id:(X,h)\to (X,k)$ induces a natural bijective map $\varphi:\partial_h^*X\to \partial_k^*X$ between the boundaries at infinity of these two Gromov hyperbolic metric spaces. On the other hand, according to Theorem \ref{zz1}, we know that $\id:(X,h)\to (X,k)$ is a rough quasi-isometry. Moreover,  Lemma \ref{g-3} guarantees that $(X,h)$ is BuSc-visual with respect to a point $w\in X$ satisfying $d(x)\leq 2d(w)$ for all $x\in X$. So appealing these facts to Theorem \ref{zz0}, it follows that $(X,k)$ is BoS-visual and thus BuSc-visual with respect to $w$, because $(X,k)$ is a length space and so is roughly geodesic. Hence we are done.
\epf

The following are our main results of this subsection that are the necessity parts of Theorems \ref{thm-2} and \ref{thm-3}, respectively. We split into two cases to investigate geometric properties of uniform metric spaces in terms of the boundedness of the spaces.

\begin{lem}\label{zz3} Let $(X,d)$ be a bounded $A$-uniform space with $A\geq1$. Then
\begin{enumerate}
\item  $(X,d)$ is locally $(\lambda,Q)$-quasiconvex for some numbers $\lambda\in(0,1)$ and $Q\geq1$;
\item  there is a point $w\in X$ such that $\diam (X)\leq C d(w)$ for some constant $C\geq 1$;
\item  $(X,k)$ is $\delta$-hyperbolic and $K$-BoS-visual with respect to $w$, for some $\delta,K\geq 0$;
\item \label{lemzz3d} the natural map $\varphi:(X_k^*,k_{w,\varepsilon})\to (\overline{X},d)$ exists and is relatively $\eta$-quasisymmetric for some homeomorphism $\eta:[0,\infty)\to [0,\infty)$.
\end{enumerate}
  The constants $\lambda,Q,C,K,\delta,\varepsilon$ and the function $\eta$ depend only on $A$.
\end{lem}

\bpf Since $X$ is $A$-uniform, it is evident that $X$ is locally $(\lambda,A)$-quasiconvex for all $\lambda\in(0,1/2)$ and the first statement holds. To verify the second one, for any number $c_0> 1$ be given, choose a point $w\in X$ with $d(x)\leq c_0 d(w)$. Then for each pair of points $x,y\in X$, since $X$ is $A$-uniform, there is an $A$-uniform arc $\gamma$ joining $x$ to $y$. Let $z\in\gamma$ be such that $\ell(\gamma[x,z])=\ell(\gamma[z,y])$. By the uniformity of $\gamma$ and by our choice of $w$, we see that
$$d(x,y)\leq \ell(\gamma)= 2\ell(\gamma[x,z])\leq 2Ad(z)\leq 2Ac_0 d(w).$$
The desired conclusion (2) follows if we take $c_0=2$ and $C=4A$.

It follows from Lemmas \ref{g-0} and \ref{h-03} that $(X,k)$ is $\delta$-hyperbolic and $K$-BoS-visual with respect to $w$, where $\delta$ and $K$ depend only on $A$. Thus (3) is true.

It remains to prove $(4)$. To this end, we see from Lemma \ref{g-3} that $(X,h)$ is BuSc-visual with respect to $w$. By Theorem \ref{zz1}, the identity map $\id:(X,h)\to (X,k)$ is roughly quasi-isometric. Moreover, Lemma \ref{g-1} implies that the natural map
$\varphi_h:(X_h^*,h_{w,\varepsilon'})\to (\overline{X},d)$
exists and is relatively $\eta_1$-quasisymmetric with $\eta_1(t)=48At^{\varepsilon'}$ for $0<\varepsilon'\leq 1/2$. On the other hand, according to Lemma \ref{zz-08}, there is a natural bijective map
$\varphi_k:X_k^*\to \overline{X}.$
It follows that the identity map $\id:(X,h)\to (X,k)$ induces a natural bijective map
$(\varphi_k|_{\partial_k^*X})^{-1}\circ \varphi_h|_{\partial_h^*X}: \partial_h^*X\to \partial_k^*X.$

Now, Lemma \ref{w-0} ensures that the identity map $\id:(X,h)\to (X,k)$ induces a map
$\id^*:(X_h^*,h_{w,\varepsilon'})\to  (X_k^*,k_{w,\varepsilon})$
which is relatively $\eta_2$-quasisymmetric for some homeomorphism $\eta_2:[0,\infty)\to [0,\infty)$, where $\varepsilon>0$ satisfies $e^{\varepsilon\delta}\leq 2$, because $(X,k)$ is roughly geodesic.
Clearly, $\varphi$ defined by $\varphi=\varphi_h\circ (\id^*)^{-1}:(X_k^*,k_{w,\varepsilon})\to (\overline{X},d)$ is relatively $\eta$-quasisymmetric with $\eta$ depending only on $\eta_1$ and $\eta_2$. Hence (\ref{lemzz3d}) holds.
\epf

\begin{lem}\label{zz4} Let $(X,d)$ be an unbounded, locally compact and $A$-uniform metric space for some constant $A\geq1$. Then
\begin{enumerate}
\item $(X,d)$ is locally $(\lambda,Q)$-quasiconvex for some numbers $\lambda\in(0,1)$ and $Q\geq1$;
\item $(X,k)$ is $\delta$-hyperbolic and $K$-roughly starlike with respect to $\xi\in\partial_k^*X$ for some $\delta,K\geq 0$;
\item\label{lemzz4c} the natural map $\psi:(X_k^*\setminus\{\xi\},k_{b,\varepsilon})\to (\overline{X},d)$ exists with $\psi(\xi)=\infty$ and is relatively $\eta$-quasisymmetric for some homeomorphism $\eta:[0,\infty)\to [0,\infty)$, where $b=b_{\xi,w}\in \mathcal{B}_k(\xi)$ and $w\in X$.
\end{enumerate}
 The constants $\lambda,Q,K,\delta,\varepsilon$ and the function $\eta$ depend only on $A$.
\end{lem}
\bpf We first observe that $(X,d)$ is locally $(\lambda,A)$-quasiconvex for all $\lambda\in(0,1/2)$, because $X$ is $A$-uniform. The $\delta$-hyperbolicity in the sense of Gromov follows from Lemma \ref{g-0} with constant $\delta=\delta(A)\geq 0$. Moreover, Lemma \ref{z2} ensures that there is a non-negative constant $K=K(A)$ such that $(X,k)$ is $K$-roughly starlike with respect to some point $\xi\in\partial_k^*X$.

It remains to check (\ref{lemzz4c}). On the one hand, we know from Lemma \ref{g-1} that the natural map
$$\varphi_h:(X_h^*\setminus\{\zeta\},k_{\beta,\varepsilon'})\to (\overline{X},d)$$
exists with $\varphi_h(\zeta)=\infty$ and is a bijection, which is relatively $\eta_1$-quasisymmetric with $\eta_1(t)=64t^{\varepsilon'},$ where $\zeta\in\partial_h^*X$ and the Busemann function $\beta=\beta_{\zeta,w}\in \mathcal{B}_h(\zeta)$ with $w\in X$. Since $X$ is $A$-uniform, it follows from Theorems \ref{zz1} and \ref{zz2} that the identity map $\id:(X,h)\to (X,k)$ is roughly quasi-isometric. Moreover, one observes that $(X,k)$ is geodesic and $(X,h)$ is BuSc-visual with respect to $w$ by Lemma \ref{z-2}.

On the other hand, Lemma \ref{zz-08} ensures that the natural map
$$\psi:X_k^*\to \overline{X}\cup\{\infty\}$$
exists with $\psi(\xi)=\infty$ and is a bijection. Thus, the identity map $\id:(X,h)\to (X,k)$ induces a bijective map
$$\phi=\psi^{-1}\circ \varphi_h:X^*_h\to X^*_k$$
with $\phi(\zeta)=\xi\in\partial^*_k X$. Therefore, we see from Lemma \ref{w-0} that
$$\phi:(X_h^*\setminus\{\zeta\},h_{\beta,\varepsilon'})\to (X_k^*\setminus\{\xi\},k_{b,\varepsilon})$$
is relatively $\eta_2$-quasisymmetric, where the Busemann function $b=b_{\xi,w}\in \mathcal{B}_k(\xi)$. Furthermore, we know that the natural map
$$\psi=\psi_h\circ (\phi)^{-1}:(X_k^*\setminus\{\xi\},k_{b,\varepsilon})\to (\overline{X},d)$$
is relatively $\eta$-quasisymmetric with $\eta$ depending only on $\eta_1$ and $\eta_2$. The proof of the lemma is complete.
\epf
\section{Uniformity from Gromov hyperbolicity}

We begin with a basic fact. Let $(X,d)$ be an incomplete, locally compact, locally quasiconvex and rectifiably connected space. Then the space $(X, k)$ is complete, proper and geodesic (cf. \cite[Proposition $2.8$]{BHK}).

\subsection{Bonk-Heinonen-Koskela's Uniformization}
We recall the following conformal deformations which were introduced by Bonk et al. (cf. \cite[Chapter $4$]{BHK}). Fix a base point $w\in X$, and consider the family of conformal deformations of $(X,k)$ induced by the densities $\rho_\epsilon(x)=e^{-\epsilon k(x,w)}$ $(\epsilon>0).$ For $x$, $y\in X$, let
$$d_\epsilon(x,y)=\inf\int_{\gamma} \rho_\epsilon \;ds_k,$$
where the infimum is taken over all rectifiable curves $\gamma$ in $X$ with endpoints $x$ and $y$, and $ds_k$ is the arc-length element with respect to $k$.
Then $d_\epsilon$ are metrics on $X$, and we  denote the resulting metric spaces by $X_\epsilon=(X,d_\epsilon)$.

Next, we cite certain auxiliary results that we shall need.

\begin{lem}\label{lem-1}$($\cite[Chapter $4$]{BHK}$)$
Suppose that $(X,d)$ is an incomplete, locally quasiconvex and rectifiably connected, locally compact space. If $(X,k)$ is $\delta$-hyperbolic and $K$-roughly starlike with respect to $w\in X$ for some $\delta,K\geq 0$, then there are numbers $A\geq 1, 0<\epsilon,c_0<1$ depending only on $\delta, K$ such that
\begin{enumerate}
  \item the space $X_\epsilon$ is $A$-uniform and has diameter at most $2/\epsilon$,
  \item for all $x, y\in X$, $c_0\epsilon k(x,y)\leq k_\epsilon(x,y)\leq e \epsilon k(x,y)$, where $k_\epsilon$ is the quasihyperbolic metric of $(X,d_\epsilon)$.
\end{enumerate}
\end{lem}

The following lemma is a generalization of \cite[Theorem $6.14$]{Vai-4}.

\begin{lem}\label{lem-2} A relatively quasisymmetric map between connected metric spaces is  relatively power quasisymmetric.
\end{lem}

\bpf Suppose that $f:(\overline{X},d)\to (\overline{Y},d')$ is relatively $\eta$-quasisymmetric and $X$ is connected. Without loss of generality, we may assume $\eta(2)\geq 2$. Since $\eta:[0,\infty)\rightarrow [0,\infty)$ is an increasing homeomorphism, there is a $\kappa\in(0,1/2)$ such that $\eta(\kappa)\leq 1/2$.

By \cite[Lemma $6.17$]{Vai-4}, we only need to prove that for all $x\in\partial X$ and $y,z\in \overline{X}$ with $d(x,y)=td(x,z)$, there are constants $C\geq 1$ and $\alpha\geq 1$ such that
\be\label{z-1.1} d'(f(x),f(y))\leq C(t^\alpha\vee t^{1/\alpha}) d'(f(x),f(z)).\ee
Towards this end, we consider two cases.

For the first one, we assume that $t\geq 1$. Note that $X$ is connected and so is the metric completion $\overline{X}$ of $X$. There is a sequence of points $y=x_0,...,x_m$ in $\overline{X}$ such that
$d(x_i,x)=2^{-i}d(y,x)$ and $2^{-m-1}d(y,x)< d(z,x)\leq 2^{-m}d(y,x).$
Put $x_{m+1}=b$. Because $f$ is relatively $\eta$-quasisymmetric, for all $1\leq i\leq m+1$, we have $$d'(f(x_{i-1}),f(x))\leq \eta(2)d'(f(x_{i}),f(x)),$$
which leads to
$$d'(f(x),f(y))=d'(f(x),f(x_0))\leq \eta(2)d'(f(x),f(x_m))\leq \eta(2)^{m+1}d'(f(x),f(z)).$$
This, together with the fact that $t\geq 2^m$, guarantees that
$$d'(f(x),f(y))\leq \eta(2)t^{\log_2\eta(2)}d'(f(x),f(z)).$$
Hence, (\ref{z-1.1}) holds in this case.

Next, we assume that $0<t<1$. Choose another sequence of points $z=y_0,\ldots,y_n$ in $\overline{X}$ such that
$d(y_j,x)=\kappa^{j}d(z,x)$ and
$$\kappa^{n+1}d(z,x)< d(y,x)\leq \kappa^{n}d(z,x)=d(y_n,x).$$
Therefore, the relative quasisymmetry of $f$ ensures that
$$d'(f(y),f(x))\leq \eta(1)d'(f(y_{n}),f(x))$$
and
$$d'(f(y_j),f(x))\leq \eta(\kappa)d'(f(y_{j-1}),f(x))\leq d'(f(y_{j-1}),f(x))/2$$
for all $1\leq j\leq n$. These two inequalities imply that
$$d'(f(y),f(x))\leq \eta(1)d'(f(y_{n}),f(x))\leq \eta(1)2^{-n}d'(f(z),f(x)).$$
Furthermore, since $t\geq \kappa^{n+1}$, we get
$$d'(f(y),f(x))\leq 2\eta(1)t^{\log_{1/\kappa} 2}d'(f(z),f(x)).$$
Hence, by taking $C=2\eta(2)$ and $\alpha=\max\{\log_2\eta(2),\log_2 1/\kappa \},$
we obtain (\ref{z-1.1}). This completes the proof.
\epf

Now we show that every relatively power quasisymmetric map is roughly quasi-isometric with respect to the hyperbolization metric introduced by Ibragimov.

\begin{lem}\label{lem-3} Let $(X,d)$ and $(Y,d')$ be two incomplete and locally compact metric spaces. If $f:\overline{X}\to \overline{Y}$ is a relatively power $\eta$-quasisymmetric bijective map with $\eta(t)=C(t^\alpha\vee t^{1/\alpha})$ for some constants $C,\alpha\geq 1$, then both $f:(X,h)\to (Y,h')$ and $f:(X,j)\to (Y,j')$  are roughly quasi-isometric.
\end{lem}
\bpf By Theorem \ref{zz2}, we only need to show that $f:(X,j)\to (Y,j')$ is rough quasi-isometric. To this end, for any two points $x,y\in X$, we may assume without loss of generality that $d'(f(x))\leq d'(f(y))$. Take a point $u\in \partial X$ with $d'(f(x),f(u))\leq 2d'(f(x))$. Set
$$d(x,y)=t(d(x)\wedge d(y)).$$
By the relative quasisymmetry of $f$, it follows that
\beq\nonumber \frac{d'(f(x),f(y))}{d'(f(x))} &\leq& 2\frac{d'(f(x),f(u))+d'(f(u),f(y))}{d'(f(x),f(u))}=2+2 \frac{d'(f(u),f(y))}{d'(f(x),f(u))}
\\ \nonumber  &\leq& 2+2\eta\Big(\frac{d(u,y)}{d(u,x)}\Big)\leq 2+2\eta\Big(\frac{d(u,x)+d(x,y)}{d(u,x)}\Big)
\\ \nonumber  &\leq& 2+2\eta\Big(1+\frac{d(x,y)}{d(x)\wedge d(y)}\Big)
\\ \nonumber  &=& 2+2C(1+t)^\alpha
\eeq
and therefore,
\beq\nonumber j'(f(x),f(y)) &\leq& \log[3+2C(1+t)^\alpha]
\\ \nonumber  &\leq& \alpha \log(1+t)+\log(5C)
\\ \nonumber  &=& \alpha j(x,y)+ \log(5C).
\eeq

Since $f$ is a bijection, its inverse $f^{-1}$ is relatively $\eta'$-quasisymmetric with $\eta'(t)=\eta^{-1}(t^{-1})^{-1}$ and thus is a relatively power quasisymmetric as well. Using a similar argument as above, we obtain the other direction. This proves that $f:(X,j)\to (Y,j')$ is roughly quasi-isometric. The lemma follows.
\epf

\subsection{Proof of Theorem \ref{thm-2}.}
Since the necessity follows from Lemma \ref{zz3}, it suffices to show the sufficiency.

For the convenience of the reader, we recall the assumptions. Suppose that $(X,d)$ is locally $(\lambda,Q)$-quasiconvex and  there is a point $w\in X$ such
that $\diam (X)\leq C \dist(w,\partial X)$. The space $(X,k)$ is  $\delta$-hyperbolic and $K$-BoS-visual with respect to $w$, and
the natural map $\varphi:(X_k^*,k_{w,\varepsilon})\to (\overline{X},d)$ exists and is relatively $\eta$-quasisymmetric.

By the assumption, we know that $X$ is a bounded space. Since  $(X,k)$ is a proper geodesic metric space, one finds from Lemma \ref{g-5} that it is $K'$-roughly starlike with respect to $w$; where $K'$ depends only on $K$ and $\delta$. Let $\epsilon>0$,
and consider the conformal deformation of $(X,k)$ and denote the conformal density by $\rho_\epsilon(x)=e^{-\epsilon k(x,w)}$ for all $x\in X.$
By Lemma \ref{lem-1}, it follows that there is a metric $d'=d_\epsilon$ induced by $\rho_\epsilon$ such that $(X,d')$ is $B$-uniform with $\diam (X,d')\leq 2/\epsilon$, and there are positive numbers $\epsilon$ and $L$ depending only on $\delta$ and $K'$ such that
\bcl\label{l-1} The identity map $\id:(X,k)\to (X,k')$ is $L$-bilipschitz, where $k'$ is the quasihyperbolic metric of $(X,d')$.
\ecl

Moreover, we prove that
\bcl\label{l-2}
 The identity map $\id:(\overline{X},d)\to (\overline{X},d')$ is relatively $\eta_0$-quasisymmetric
for some homeomorphism $\eta_0:[0,\infty)\to [0,\infty)$.
\ecl

This can be seen as follows. By Claim \ref{l-1} and Lemma \ref{zz}, the identity map $\id:(X,k)\to (X,k')$ induces  a relatively $\eta_1$-quasisymmetric map
$$\phi:(X_k^*,k_{w,\varepsilon})\to (X_{k'}^*,k'_{w,\varepsilon'}).$$
By \cite[$(4.17)$]{BHK}, we have
$$d'(w)\geq \frac{\rho_\epsilon(w)}{\epsilon e}\geq \frac{1}{2e}\diam(X,d').$$
Since $(X,d')$ is $B$-uniform,  Lemma \ref{zz3} ensures that the natural map
$$\psi:(X_{k'}^*,k'_{w,\varepsilon'})\to (\overline{X},d')$$
exists and is relatively $\eta_2$-quasisymmetric. Therefore, using the above facts, we obtain that the identity map
$$\id=\psi\circ \phi\circ \varphi^{-1}:(\overline{X},d)\to (\overline{X},d')$$
is a relatively $\eta_0$-quasisymmetric map with $\eta_0=\eta_2\circ \eta_1\circ \eta^{-1}$. Hence, Claim \ref{l-2} is valid.

Now, it follows from Claim \ref{l-2}, Lemmas \ref{lem-2} and \ref{lem-3} that the identity map $\id:(X,j)\to (X,j')$ is roughly quasi-isometric. Since $(X,d')$ is $B$-uniform, we see from Theorem \ref{zz1} that the identity map $\id:(X,k')\to (X,j')$ is roughly quasi-isometric as well. Thus by Claim \ref{l-1}, we find that the identity map $\id:(X,k)\to (X,j)$ is also roughly quasi-isometric. Again by Theorem \ref{zz1}, it follows that $(X,d)$ is uniform. The proof of Theorem \ref{thm-2} is complete.
\qed

\subsection{Uniformization II}
Before proving Theorem \ref{thm-3}, we recall some facts concerning the theory of unbounded uniformizing Gromov hyperbolic spaces (cf. \cite{Z}).

Let $(X,d)$ be an incomplete, rectifiably connected, locally compact and locally quasiconvex metric space, and let $k$ be the quasihyperbolic metric of $X$. Fix $w\in X$, $\xi\in\partial_\infty X$, and  take a Busemann function $b=b_{\xi,w}\in \mathcal{B}_k(\xi)$. Consider the family of conformal deformations of $(X, k)$ by the densities $\varrho_\epsilon(x)=e^{-\epsilon b(x)}$ for $\epsilon>0.$
We denote the resulting metric spaces by $(X,\chi_\epsilon)$. Thus $\chi_\epsilon$ is a metric defined by
$$\chi_\epsilon(x,y)=\inf \int_\gamma \varrho_\epsilon \;ds_k,$$
where the infimum is taken over all rectifiable curves in $X$ joining the points $x$ and $y$. Then we cite the following useful result.

\begin{lem}\label{lem-4} $($\cite[Theorem 1.2 and Lemma 5.24]{Z}$)$
Suppose that $(X,d)$ is an incomplete, locally quasiconvex, rectifiably connected and locally compact space. If $(X,k)$ is $\delta$-hyperbolic and $K$-roughly starlike with respect to $\xi\in \partial^*_k X$, then there are positive numbers $A, \epsilon, L$ depending only on $\delta$ and $K$ such that
\begin{enumerate}
\item  The density $\varrho_\epsilon$ induces an unbounded $A$-uniform metric space $X_\epsilon$ such that any sequence $\{x_n\}\in\xi$ in $X_\epsilon$ tends to infinity;
\item The identity map $(X, k)\to (X,k_\epsilon)$
is $L$-bilipschitz, where $k_\epsilon$ is the quasihyperbolic metric of $(X,\chi_\epsilon)$.
 \end{enumerate}
\end{lem}

\subsection{Proof of Theorem \ref{thm-3}.}
Because the necessity follows from Lemma \ref{zz4}, we only need to show the sufficiency. {First, we recall the assumptions.} Suppose that $(X,d)$ is a locally compact, rectifiably connected, locally $(\lambda,Q)$-quasiconvex, and incomplete metric space, and $(X,k)$ is  $\delta$-hyperbolic and $K$-roughly starlike with respect to $\xi\in\partial_k^*X$. Moreover, the natural map
$$\psi:(X_k^*\setminus\{\xi\},k_{b,\varepsilon})\to (\overline{X},d)$$
exists with $\psi(\xi)=\infty$ and is relatively $\eta$-quasisymmetric, where the Busemann function $b=b_{\xi,w}\in \mathcal{B}_k(\xi)$ and $w\in X$.

Let $\epsilon>0$, and consider the conformal deformation of $(X,k)$ by the density $\varrho_\epsilon(x)=e^{-\epsilon b(x)}$ for all $x\in X,$
where $b=b_{\xi,w}\in \mathcal{B}_k(\xi)$ and $w\in X$. By Lemma \ref{lem-4}, we find that there is a metric $\widetilde{d}=\chi_\epsilon$ induced by $\varrho_\epsilon$ such that $(X,\widetilde{d})$ is an unbounded $C$-uniform space with any sequence $\{x_n\}\in\xi$ in $(X,\widetilde{d})$ tending to infinity. Moreover, the identity map
\be\label{l-3} \id:(X,k)\to (X,\widetilde{k})\;
\ee
is $L$-bilipschitz for some positive numbers $\epsilon$ and $L$ depending only on $\delta$ and $K$, where $\widetilde{k}$ is the quasihyperbolic metric of $(X,\widetilde{d})$.

We next prove that the identity map
\be\label{l-4} \id:(\overline{X},d)\to (\overline{X},\widetilde{d})\;
\ee
is relatively $\eta_0$-quasisymmetric for some homeomorphism $\eta_0:[0,\infty)\to [0,\infty)$.

To this end, it follows from the above observation about $(\ref{l-3})$ and Lemma \ref{zz} that there is a relatively $\eta_1$-quasisymmetric map
$$\phi:(X_k^*\setminus\{\xi\},k_{b,\varepsilon})\to (X_{\widetilde{k}}^*\setminus\{\zeta\},\widetilde{k}_{\beta,\varepsilon'})$$
induced by the identity map $\id:(X,k)\to (X,\widetilde{k})$. Here $\beta$ is a Busemann function with $\beta=\beta_{\xi,w}\in \mathcal{B}_{\widetilde{k}}(\zeta)$ with $\phi(\xi)=\zeta$. Since $(X,\widetilde{d})$ is $C$-uniform,  Lemma \ref{zz4} guarantees that the natural map
$$\widetilde{\psi}: (X_{\widetilde{k}}^*\setminus\{\zeta\},\widetilde{k}_{\beta,\varepsilon'}) \to (\overline{X},\widetilde{d})$$
exists with $\psi(\zeta)=\infty$ and is relatively $\eta_2$-quasisymmetric. Therefore, using the above facts, we know that the identity map
$$\id=\widetilde{\psi}\circ \phi\circ \psi^{-1}:(\overline{X},d)\to (\overline{X},\widetilde{d})$$
is relatively $\eta_0$-quasisymmetric with $\eta_0=\eta_2\circ \eta_1\circ \eta^{-1}$. Hence, the identity map given by \eqref{l-4} is relatively $\eta_0$-quasisymmetric.

Finally, in view of the above facts, a similar argument as in the last part of the proof of Theorem \ref{thm-2} ensures that $(X,d)$ is uniform. This completes the proof.
\qed

\section{Uniform spaces and Teichm\"uller's displacement property}

\subsection{Teichm\"uller's displacement property and Gromov hyperbolicity imply uniformity}
For the convenience of the reader, we recall the following definition from Subsection \ref{sec-t}. Let $(X,d)$ be a rectifiably connected, locally compact and incomplete metric space. We define
\beq\nonumber \mathcal{T}_{M,\mu}(\overline{X})
&=& \Big\{{X} \xrightarrow{f} {X}:(X,h)\xrightarrow{f}(X,k)\; \mbox{is an $(M,\mu)$-rough quasi-isometry},
\\ \nonumber & & ~~~~ \mbox{ the continuous extension of $f$ on the boundary}\; f|_{\partial X}=\id_{\partial X} \Big\}.
\eeq

We first consider the bounded case and our result is as follows.

\begin{lem}\label{zz-05} Let $(X,d)$ be a bounded, locally compact, locally $(\lambda,Q)$-quasiconvex, rectifiably connected and incomplete metric space. If
\begin{enumerate}
\item\label{q-2} $(X,k)$ is $\delta$-hyperbolic, and the natural map $\varphi:(\partial_k^* X,k_{w,\varepsilon})\to (\partial{X},d)$ exists and is power $\eta$-quasisymmetric;
\item\label{q-3} $(X,d)$ satisfies the $\Lambda$-displacement property with respect to the quasihyperbolic metric $k$,
\end{enumerate}
then $(X,d)$ is $A$-uniform with $A$ depending only on $\lambda,Q,\delta,\eta,\Lambda$ and $\varepsilon$.
\end{lem}

\bpf
By Theorem \ref{zz2}, we observe that the natural map
$\psi:\partial_h^*X\to \partial X$
exists and is a bijection. It follows from this fact and the first condition (\ref{q-2}) that the identity map $\id:(X,h)\to (X,k)$ induces a bijective map
$$\phi=\varphi^{-1}\circ \psi:\partial_h^*X\to \partial_k^*X.$$
By Lemma \ref{g-3}, we know that $(X,h)$ is a BuSc-visual hyperbolic space. Because $(X,k)$ is a roughly geodesic Gromov hyperbolic space, Lemma \ref{z6} ensures that there is an $(M,\mu)$-roughly quasi-isometric map
$f:(X,h)\to (X,k)$
with $f|_{\partial_h^*X}=\phi$. Thus we have
$$f|_{\partial X}=\varphi\circ \phi\circ\psi^{-1}=\id_{\partial X},$$
which implies that $f\in \mathcal{T}_{M,\mu}(\overline{X})$. By invoking the second condition $(\ref{q-3})$, it follows that there is a constant $C\geq 0$ depending only on $M$, $\mu$ and the data of $X$ such that $k(x,f(x))\leq C$ for all $x\in X$.
Therefore, we have for all $x,y\in X$
$$k(x,y)\leq k(x,f(x))+k(f(x),f(y))+k(f(y),y)\leq Mh(x,y)+\mu+2C.$$
This, together with Theorems \ref{zz1} and \ref{zz2}, shows that $(X,d)$ is uniform.
\epf

For the unbounded case, we have the following result.

\begin{lem}\label{zz-10} Let $(X,d)$ be an unbounded, locally compact, locally $(\lambda,Q)$-quasiconvex, incomplete and rectifiably connected metric space. If
\begin{enumerate}
\item\label{q-4} $(X,k)$ is $\delta$-hyperbolic, and the natural map $\psi:(\partial_k^* X\setminus\{\xi\},k_{b,\varepsilon})\to (\partial{X},d)$ exists and is power $\eta$-quasisymmetric, where $\psi(\xi)=\infty$ and $b=b_{\xi,w}\in \mathcal{B}_k(\xi)$ is a Busemann function with $w\in X$;
\item\label{q-5} $(X,d)$ satisfies the $\Lambda$-displacement property with respect to the quasihyperbolic metric $k$,
\end{enumerate}
then $(X,d)$ is $A$-uniform with $A$ depending only on $\lambda,Q,\delta,\eta,\Lambda$ and $\varepsilon$.
\end{lem}
\bpf Since $X$ is unbounded, by Lemma \ref{g-1}, we see that the natural map
$$\psi_h:(\partial_h^*X\setminus \{\zeta\},h_{\beta,\varepsilon'})\to(\partial X,d)$$
exists and is power quasisymmetric, where $\psi_h(\xi')=\infty=\psi(\xi)$ and $\beta=\beta_{\zeta,w}\in \mathcal{B}_h(\zeta)$. Using this fact and the first condition $(\ref{q-4})$, we obtain a power quasisymmetric map
$$\phi=\psi^{-1}\circ \psi_h:(\partial_h^*X\setminus \{\xi'\},h_{\beta,\varepsilon'})\to(\partial_k^* X\setminus\{\xi\},k_{b,\varepsilon}).$$

On the other hand, Lemma \ref{g-3} ensures that $(X,h)$ is a BuSc-visual hyperbolic space. Since $(X,k)$ is a roughly geodesic Gromov hyperbolic space, we observe from Lemma \ref{z6} that there is an $(M,\mu)$-roughly quasi-isometric map
$f:(X,h)\to (X,k)$
with $f^*|_{\partial_h^*X}=\phi$. Thus
$$f|_{\partial X}=\psi\circ \phi\circ \psi_h^{-1}=\id_{\partial X}.$$
These two facts imply that $f\in \mathcal{T}_{M,\mu}(\overline{X})$.  By the second condition (\ref{q-5}), it follows that there is a constant $C\geq 0$ depending only on $M,\mu$ and the data of $X$ such that  $k(x,f(x))\leq C$ for all $x\in X$.
Therefore, for all $x\in X$,  we obtain
$$k(x,y)\leq k(x,f(x))+k(f(x),f(y))+k(f(y),y)\leq Mh(x,y)+\mu+2C.$$
This inequality together with Theorems \ref{zz1} and \ref{zz2}, asserts that $(X,d)$ is uniform and the proof is complete.
\epf

\subsection{Uniform spaces satisfy Teichm\"uller's displacement property} To get a quantitative result, we recall some notation introduced by Bonk and  Kleiner in \cite{BK1}:
$$[x,y,z,w]=\frac{d(x,z)\wedge d(y,w)}{d(x,y)\wedge d(z,w)},$$
for any points $x,y,z,w$ in a metric space $(X,d)$. They also established a relation between the cross ratio $\tau(x,y,z,w)$ and $[x,y,z,w]$ which reads as follows.

\begin{lem}\label{zz-01} $($\cite[Lemma 2.3]{BK1}$)$  For any points $x,y,z,w$ in a metric space $(X,d)$, we have
$[x,y,z,w] \leq \theta_0(\tau(x,y,z,w)),$
where $\theta_0(t)=3(t\vee \sqrt{t})$.
\end{lem}

We divide the proof into two parts to show that uniform spaces satisfy Teichm\"uller's displacement property.

\begin{lem}\label{zz-06} Suppose that $(X,d)$ is a bounded $A$-uniform locally compact metric space with $C$-uniformly perfect boundary and $\diam (X)\leq \nu\diam (\partial X)$ for some constants $A,C,\nu\geq 1$. Then $(X,d)$ satisfies the displacement property with respect to the quasihyperbolic metric $k$, quantitatively.
\end{lem}
\bpf Because $X$ is $A$-uniform, it follows from Theorems \ref{zz1} and \ref{zz2} that the identity map $\id:(X,h)\to (X,k)$ is roughly quasi-isometric. For a given map $f\in \mathcal{T}_{M,\mu}(\overline{X})$, we see that the mapping $f:(X,k)\to (X,k)$ is also roughly quasi-isometric.

To check the displacement property of $X$ with respect to the quasihyperbolic metric $k$, by Theorem \ref{zz1}, it suffices to find a function $\Lambda$ depending on $A,C,\nu$ and $M,\mu$ such that
\be\label{z-1.0} j(x,f(x))\leq \Lambda \ \mbox{ for all $x\in X$.}\ee

To establish (\ref{z-1.0}), we need some preparation. First, we prove that

\bcl\label{z-1.2} $f:\overline{X}\to \overline{X}$ is relatively $\theta_1$-quasim\"{o}bius with $f|_{\partial X}=\id_{\partial X}$. \ecl

We only need to show the first assertion, because $f\in \mathcal{T}_{M,\mu}(\overline{X})$ implies that $f|_{\partial X}=\id_{\partial X}$. To this end, since $X$ is a bounded $A$-uniform metric space, Lemma \ref{zz-08} guarantees that the natural map $\varphi:\partial_k^*X\to \partial X$ exists and is a bijection. By \cite[Proposition $6.3$]{BS}, it follows that $f:(X,k)\to (X,k)$ induces a boundary map $\partial^*f:\partial_k^*X\to \partial_k^*X$, which is an injection. Since $f|_{\partial X}=\id_{\partial X}$, we see that
$$\partial^*f=\varphi^{-1}\circ f|_{\partial X} \circ \varphi=\id_{\partial^*_k X}:\partial_k^*X\to\partial_k^*X$$
 is a bijection.

Now, choose a point $w\in X$ with $d(x)\leq 2d(w)$ for all $x\in X$. We observe from Lemma \ref{zz} that $f$ induces a relatively quasisymmetric map
$f^*:(X_k^*,k_{w,\varepsilon})\to (X_k^*,k_{f(w),\varepsilon'}),$
with $f^*|_{\partial^*_k X}=\id_{\partial^*_k X}$. Moreover, a similar argument as that of  Lemma \ref{z-1} ensures that the identity map $(X_k^*,k_{w,\varepsilon})\to (X_k^*,k_{f(w),\varepsilon'})$ is relatively quasim\"{o}bius, we find that
$f^*:(X_k^*,k_{w,\varepsilon})\to (X_k^*,k_{w,\varepsilon'})$
is also relatively quasim\"{o}bius, because the compositions of relatively quasim\"{o}bius maps and relatively quasisymmetric maps are relatively quasim\"{o}bius.

On the other hand, since $X$ is a bounded $A$-uniform metric space, Lemma \ref{zz3} ensures that the natural map
$\psi:(X_k^*, k_{w,\varepsilon})\to (\overline{X},d)$
exists and is relatively quasisymmetric. Therefore, we find that
$f=\psi\circ f^*\circ \psi^{-1}:(\overline{X},d)\to (\overline{X},d)$
is relatively quasim\"{o}bius. Thus Claim \ref{z-1.2} holds.

\bcl\label{z-1.3} $f:\overline{X}\to \overline{X}$ is relatively $\theta$-quasisymmetric with $f|_{\partial X}=\id_{\partial X}$. \ecl

Take two points $z_1,z_2\in \partial X$ with $d(z_1,z_2)\geq \diam (\partial X)/2$. Since $\partial X$ is $C$-uniformly perfect, there is a  $z_3\in \partial X$ such that
$$(2C)^{-1}d(z_1,z_2)\leq d(z_1,z_3)\leq 2^{-1}d(z_1,z_2).$$
The triangle inequality yields that
$$d(z_2,z_3)\geq 2^{-1}d(z_1,z_2)\geq 4^{-1}\diam (\partial X).$$
Therefore, we observe that there are three points $z_i$, $i=1,2,3$ in $\partial X$, such that
\be\label{z-3} d(z_i,z_j)\geq \frac{1}{4C}\diam (\partial X) \geq \frac{1}{4C\nu}\diam (X) \ \mbox{ for all $1\leq i<j\leq 3$.} \ee

For any three distinct points $x,y,z$ in $\overline{X}$ with $x\in \partial X$, it follows from the three-point condition $(\ref{z-3})$ that, for all $u\in \overline{X}$, there are $i\neq j\in\{1,2,3\}$ such that
$$d(u,z_i)\wedge d(u,z_j)\geq \frac{1}{8C\nu}\diam (X).$$
This ensures that there exists a point $z_i\in \{z_1,z_2,z_3\}$ such that
$$d(y,z_i)\geq \frac{1}{8C\nu}\diam (X) \; \mbox{ and }\; d(f(z),z_i)\geq \frac{1}{8C\nu}\diam (X).$$
Because $f|_{\partial X}=\id_{\partial X}$, we have
$$d(y,z_i)\wedge d(z,x) \geq \frac{1}{8C\nu} d(z,x) \; \mbox{and}\; d(f(z),z_i) \wedge d(x,f(y))\geq \frac{1}{8C\nu} d(x,f(y)),$$
which guarantee that
\be\label{z-4} [x,z,y,z_i]\leq 8C\nu \frac{d(x,y)}{d(x,z)} \;\mbox{and}\; [f(x),f(z),f(y),f(z_i)] \geq \frac{d(f(x),f(y))}{8C\nu d(f(x),f(z))}. \ee

On the other hand, Claim \ref{z-1.2} implies that $f$ is relatively $\theta_1$-quasim\"{o}bius. Moreover, since $x,z_i\in \partial X$, it follows from Lemma \ref{zz-01} that
$$[f(x),f(z),f(y),f(z_i)] \leq \theta_2([x,z,y,z_i]),$$
where $\theta_2(t)=\theta_0\circ \theta_1\circ (\frac{1}{\theta_0^{-1}}(1/t))$. Now by (\ref{z-4}), we obtain
\beq\nonumber
\frac{d(f(x),f(y))}{d(f(x),f(z))} &\leq& 8C\nu [f(x),f(z),f(y),f(z_i)]
\\ \nonumber&\leq& 8C\nu \theta_2[(x,z,y,z_i)]
\\ \nonumber&\leq&8C\nu \theta_2 \bigg(8C\nu \frac{d(x,y)}{d(x,z)}\bigg).
\eeq
By taking $\theta(t)=8C\nu \theta_2 (8C\nu t)$, it follows that $f$ is relatively $\theta$-quasisymmetric.

\bcl\label{z-1.4}
For each $x\in X$ and $x_0\in\partial X$ with $d(x,x_0)\leq 2d(x)$, there is a point $x_1\in\partial X$ such that $(2C\nu)^{-1}d(x) \leq d(x_0,x_1) \leq 5d(x).$
\ecl
We consider two cases. If $\partial X\subset \overline{B}(x,3d(x))$, then take a point $x_1\in\partial X$ with $d(x_0,x_1)\geq 2^{-1}\diam (\partial X)$.
Then we obtain
$$d(x_0,x_1)\geq \frac{\diam(\partial X)}{2\diam (X)}d(x)\geq \frac{1}{2\nu} d(x)$$
and
$$d(x_0,x_1)\leq d(x_0,x)+d(x,x_1)\leq 5d(x),
$$
as required. For the other case, that is, $\partial X\setminus \overline{B}(x,3d(x))\not=\emptyset$, by the triangle inequality, we have $\overline{B}(x_0,d(x))\subset \overline{B}(x,3d(x)),$
and thus,
$\partial X\setminus \overline{B}(x_0,d(x))\not=\emptyset.$
Moreover, since $\partial X$ is $C$-uniformly perfect, there is a point $x_1\in\partial X$ such that
$$\frac{d(x)}{2C}\leq d(x_0,x_1)\leq d(x)/2,
$$
and Claim \ref{z-1.4} follows.

\bcl\label{z-1.5}
For any $x\in X$ and $x_0\in\partial X$ with $d(x,x_0)\leq 2d(x)$, there is a constant $L\geq 1$
such that
$$d(x)/L\leq d(f(x))\leq d(f(x),f(x_0))\leq Ld(x).$$
\ecl
For the one direction, we see from Claim \ref{z-1.4} that there is a point $x_1\in\partial X$ such that
$\frac{1}{2C\nu}d(x) \leq d(x_0,x_1) \leq 5d(x).$
Claim \ref{z-1.3} ensures that $f$ is relatively $\theta$-quasisymmetric with $f|_{\partial X}=\id_{\partial X}$. Moreover, for the triplet of points $x,x_0,x_1\in \overline{X}$, we get
$$ \frac{d(f(x_0),f(x))}{d(f(x_0),f(x_1))}\leq \theta\bigg(\frac{d(x_0,x)}{d(x_0,x_1)}\bigg)\leq \theta(4C\nu),$$
and thus,
\be\label{z-0.1} d(f(x))\leq d(f(x_0),f(x))\leq \theta(4C\nu)d(x_0,x_1)\leq 5\theta(4C\nu)d(x).\ee

For the other direction, we discuss by a symmetric argument. Choose a point $x_2\in \partial X$ with $d(f(x),f(x_2))\leq 2d(f(x))$. Again by Claim \ref{z-1.4}, it follows that there is a point $x_3\in\partial X$ such that
$$\frac{1}{2C\nu}d(f(x)) \leq d(f(x_2),f(x_3)) \leq 5d(f(x)).$$
Since $f$ is relatively $\theta$-quasisymmetric with $f|_{\partial X}=\id_{\partial X}$, for the triplet of points $x,x_2,x_3\in \overline{X}$, we have
$$\frac{d(f(x_2),f(x_3))}{d(f(x_2),f(x))}\leq \theta\bigg(\frac{d(x_2,x_3)}{d(x_2,x)}\bigg).$$
Thus a direct computation gives that
$$\frac{d(x_2,x)}{d(x_2,x_3)}\leq \widetilde{\theta}\bigg(\frac{d(f(x_2),f(x))}{d(f(x_2),f(x_3))}\bigg)\leq \widetilde{\theta}(4C\nu),$$
where $\widetilde{\theta}(t)=\theta^{-1}(t^{-1})^{-1}$ for $t>0$. This implies that
\be\label{z-0.2} d(x)\leq d(x_2,x)\leq \widetilde{\theta}(4C\nu)d(x_2,x_3)\leq 5\widetilde{\theta}(4C\nu)d(f(x)).\ee
Therefore, by (\ref{z-0.1}) and (\ref{z-0.2}), Claim \ref{z-1.4}  follows with the choice of $L=5\theta(4C\nu)\vee (5\widetilde{\theta}(4C\nu)).$

Finally, we are ready to prove the inequality \eqref{z-1.0}.
Choose a point $x_0\in\partial X$ with $d(x,x_0)\leq 2d(x)$. Because $d(x,f(x))\leq d(x,x_0)+d(f(x),f(x_0))\leq 2L d(x)$, it follows from Claim \ref{z-1.5} that
$$j(x,f(x))=\log\bigg(1+\frac{d(x,f(x))}{d(x)\wedge d(f(x))}\bigg)\leq \log(1+2L^2)=\Lambda,$$
as desired. The proof of the lemma is complete.
\epf

For the unbounded case, the following result may be proved along the same lines of proof of Lemma \ref{zz-06}. For the sake of completeness, we include the details.

\begin{lem}\label{zz-09} Let $A,C\geq 1$. Suppose that $(X,d)$ is an unbounded $A$-uniform locally compact metric space, and that  $\partial X\cup\{\infty\}$ is a $C$-uniformly perfect set. Then $(X,d)$ satisfies the displacement property with respect to the quasihyperbolic metric $k$, quantitatively.
\end{lem}

\bpf
Since $X$ is $A$-uniform, we see from Theorems \ref{zz1} and \ref{zz2} that the identity map $\id:(X,h)\to (X,k)$ is roughly quasi-isometric. For a given map $f\in \mathcal{T}_{M,\mu}(\overline{X}),$ it follows that the mapping $f:(X,k)\to (X,k)$ is also roughly quasi-isometric.

To demonstrate the displacement property of $X$ with respect to the quasihyperbolic metric $k$, by Theorem \ref{zz1}, we only need to show that
\be\label{z-2.1} j(x,f(x))\leq \Lambda \ \mbox{ for all $x\in X$,}\ee
where $\Lambda$ is a constant depending only on $A$ and $C$.

\bcl\label{z-5} $f:\overline{X}\to \overline{X}$ is relatively $\eta$-quasisymmetric with $f|_{\partial X}=\id_{\partial X}$.\ecl

Since $X$ is unbounded  and $A$-uniform,  Lemma \ref{zz-08} implies that the natural map
$$\varphi:\partial_k^*X\to \partial X\cup \{\infty\}$$
exists and is a bijection with $\varphi(\xi)=\infty$ for some $\xi\in\partial_k^*X$. Then by \cite[Proposition $6.3$]{BS}, we observe that $f:(X,k)\to (X,k)$ induces an injective boundary map
$\partial^*f:\partial_k^*X\to \partial_k^*X.$
Since $f|_{\partial X}=\id_{\partial X}$, we get a bijection
$$\partial^*f=\varphi^{-1}\circ f|_{\partial X} \circ \varphi=\id_{\partial^*_k X}:\partial_k^*X\to\partial_k^*X.$$
Since $(X,k)$ is a geodesic Gromov hyperbolic metric space and  $\partial^*f=\id_{\partial^*_k X}$, Lemma \ref{zz} guarantees that
$$f_*:(X_k^*\setminus\{\xi\},\sigma_{b,\varepsilon}) \to (X_k^*\setminus\{\xi\},\sigma_{b',\varepsilon'})$$
is relatively quasisymmetric, where $b=b_{\xi,w}, b'=b'_{\xi,w'}\in \mathcal{B}_k(\xi)$ and $w'=f(w)$.

On the other hand, it follows from Theorem \ref{thm-3} that the natural map $$\psi_1:(X_k^*\setminus\{\xi\},k_{b,\varepsilon})\to (\overline{X},d)$$
exists and is a relatively $\eta_1$-quasisymmetric homeomorphism for some function $\eta_1:[0,\infty)\to [0,\infty)$. Similarly, the natural map
$$\psi_2:(X_k^*\setminus\{\xi\},k_{b',\varepsilon'})\to (\overline{X},d)$$
is a relatively $\eta_2$-quasisymmetric homeomorphism for some  $\eta_2:[0,\infty)\to [0,\infty)$ with $\psi_1(x)=\psi_2(x)=x$ for all $x\in X$. Thus we see that
$$f=\psi_2\circ f_*\circ \psi_1^{-1}: (\overline{X},d)\to (\overline{X},d)$$
is a relatively quasisymmetric map with $f|_{\partial X}=\id_{\partial X}$, because the compositions of relatively quasisymmetric maps are also relatively quasisymmetric. This proves Claim \ref{z-5}.

Let $x\in X$. Take a point $x_0\in\partial X$ with $d(x,x_0)=d(x)$. Note that the one-point extended boundary $\partial X\cup\{\infty\}$ is $C$-uniformly perfect associated with a sphericalized metric which was introduced by Bonk and Kleiner \cite{BK}. Note that the space equipped with the original metric and the sphericalized metric are quasim\"obius to each other under the identity map, and
the uniform perfectness is preserved under quasim\"obius maps in a quantitative way. With no loss of generality we may assume that $(\partial X,d)$ is also $C$-uniformly perfect. Thus there is a point $x_1\in\partial X$ such that
$(2C)^{-1}d(x)\leq d(x_0,x_1)\leq 2^{-1}d(x).$  Now, we proceed to prove the next claim.

\bcl\label{z-6} There is a constant $L\geq 1$ such that
$$d(x)/L\leq d(f(x))\leq d(f(x),f(x_0))\leq L d(x).$$ \ecl

By Claim \ref{z-5}, we know that $f$ is relatively $\eta$-quasisymmetric and $f|_{\partial X}=\id_{\partial X}$. Thus for the triplet of points $x,x_0,x_1\in \overline{X}$, we compute
$$\frac{d(f(x_0),f(x))}{d(f(x_0),f(x_1))}\leq \eta\bigg(\frac{d(x_0,x)}{d(x_0,x_1)}\bigg)\leq \eta(2C),$$
which implies
$$d(f(x))\leq d(f(x_0),f(x))\leq \eta(2C)d(x_0,x_1)\leq 2\eta(2C)d(x).$$
By symmetry, a similar argument as in the proof of Claim \ref{z-1.5} gives the estimate $d(x)\leq Ld(f(x))$ for some constant $L>0$. Hence Claim \ref{z-6} holds.

Now we check that $j(x,f(x))\leq \Lambda$ for some positive number $\Lambda$. Indeed, it follows from Claim \ref{z-6} that
$$j(x,f(x))=\log\bigg(1+\frac{d(x,f(x))}{d(x)\wedge d(f(x))}\bigg)\leq \log(1+2L^2)=\Lambda,$$
because $d(x,f(x))\leq d(x,x_0)+d(f(x),f(x_0))\leq 2L d(x)$.
\epf

\subsection{Proofs of Theorems \ref{thm-4} and \ref{thm-5}}
Theorem \ref{thm-4} follows from Lemmas \ref{zz-05} and \ref{zz-06},  where as the proof of Theorem \ref{thm-5} follows from Lemmas \ref{zz-10} and \ref{zz-09}.
\qed


\end{document}